\documentclass[11pt]{article}
\usepackage[left=1in,right=1in,top=1in,bottom=1in]{geometry}
\usepackage{times}
\usepackage{expl3}
\usepackage{cite}
\usepackage[table]{xcolor}
\usepackage{multirow}
\usepackage{stackengine} 
\usepackage{hhline}
\usepackage{lipsum}
\usepackage{titlesec}
\usepackage{wrapfig}
\usepackage{epsfig}
\usepackage{graphicx}
\usepackage{amsmath}
\usepackage{amssymb}
\usepackage{epstopdf}
\usepackage{boldline}
\usepackage{datetime} 
\usepackage{calligra}
\usepackage{url}
\usepackage{blindtext}

\newcommand{\define}{\stackrel{\mbox{\tiny def}}{=}}

\newtheorem{theorem}{Theorem}
\newtheorem{corollary}{Corollary}
\newtheorem{lemma}{Lemma}

\usepackage{mathtools}
\usepackage{epstopdf}
\usepackage{balance}
\usepackage{thmtools}
\usepackage{thm-restate}
\usepackage{hyperref}
\usepackage{cleveref}

\usepackage[ruled,vlined]{algorithm2e}
\include{pythonlisting}
\newcommand{\ostar}{\mathbin{\mathpalette\make@circled\star}}

\makeatletter
\newcommand{\removelatexerror}{\let\@latex@error\@gobble}
\makeatother
\makeatletter
\newcommand*{\rom}[1]{\expandafter\@slowromancap\romannumeral #1@}
\makeatother

\ExplSyntaxOn
\newcommand\latinabbrev[1]{
  \peek_meaning:NTF . {
    #1\@}%
  { \peek_catcode:NTF a {
      #1.\@ }%
    {#1.\@}}}
\ExplSyntaxOff

\titleclass{\subsubsubsection}{straight}[\subsubsection]

\begin{document}
\vspace{1cm}
\title{T-product Tensor Expander Chernoff Bound}\vspace{1.8cm}
\author{Shih~Yu~Chang 
\thanks{}
\thanks{Shih Yu Chang is with the Department of Applied Data Science,
San Jose State University, San Jose, CA, U. S. A. (e-mail: {\tt
shihyu.chang@sjsu.edu}).
           }}
\maketitle


\begin{abstract}
In probability theory, the Chernoff bound gives exponentially decreasing bounds on tail distributions for sums of independent random variables and such bound is applied at different fields in science and engineering. In this work, we generalize the conventional Chernoff bound from the summation of independent random variables to the summation of dependent random T-product tensors. Our main tool used at this work is majorization technique. We first apply majorizaton method to establish norm inequalitites for T-product tensors and these norm inequalities are used to derive T-product tensor expander Chernoff bound. Compared with the matrix expander Chernoff bound obtained by Garg et al., the T-product tensor expander Chernoff bound proved at this work contributes following aspects: (1) the random objects dimensions are increased from matrices (two-dimensional data array) to T-product tensors (three-dimensional data array); (2) this bound generalizes the identity map of the random objects summation to any polynomial function of the random objects summation; (3) Ky Fan norm, instead only the maximum or the minimum eigenvalues, for the function of the random T-product tensors summation is considered; (4) we remove the restriction about the summation of all mapped random objects is zero, which is required in the matrix expander Chernoff bound derivation. 
\end{abstract}


\section{Introduction}\label{sec:Introduction} 

The Chernoff bound provides the exponential decreasing inequality on tail distribution of sums of independent random variables and such bound is used extensively at different fields of science and engineering. For instance, the Chernoff bound is used to estimate approximation error of statistical machine learning algorithms~\cite{kearns1994introduction}. In communication networking system, the Chernoff bound is utilized to establish bounds for packet routing problems which are used to design congestion reduction routing protocol in sparse networks~\cite{jang1999chernoff}. It is a tighter bound than the known first- or second-moment-based tail bounds such as Markov's inequality or Chebyshev's inequality, which only yield power-law bounds on tail distribution. However, neither Markov's inequality nor Chebyshev's inequality requires that the variates are independent, which is necessary by the Chernoff bound~\cite{chernoff1981note}.

There are several directions to generalize the Chernoff bound. One major direction is to increase the dimension of random objects from random variables to random matrices. The works of Rudelson~\cite{rudelson1999random}, Ahlswede-Winter~\cite{ahlswede2002strong} and Tropp~\cite{tropp2012user} demonstrated that a similar concentration bound is also valid for matrix-valued random variables. If $\mathbf{X}_1,\mathbf{X}_2, \cdots, \mathbf{X}_n$ are independent $m \times m$ Hermitian complex random matrices with $\left\Vert  \mathbf{X}_i \right\Vert \leq 1$ for $1 \leq i \leq n$, where $\left\Vert \cdot \right\Vert$ is the spectral norm, we have following Chernoff bound for the version of $n$ i.i.d. random matrices: 
\begin{eqnarray}\label{eq:Chernoff bound iid matrices}
\mathrm{Pr} \left( \left\Vert \frac{1}{n} \sum\limits_{i=1}^{n} \mathbf{X}_i - \mathbb{E}[\mathbf{X}] \right\Vert \geq \vartheta \right) \leq m \exp( - \Omega  n \vartheta^2),
\end{eqnarray}
where $\Omega$ is a constant related to the matrix norm. This is also called ``Matrix Chernoff Bound'' and is applied to many fields, e.g., spectral graph theory, numerical linear algebra, machine learning and information theory~\cite{tropp2015introduction}. Recently, the author generalized matrix bounds to various tensors bounds under Einstein product, e.g., Chernoff, Bennett, and Bernstein inequalities associated with tensors under Einstein product in~\cite{chang2020convenient}. 

Another direction to extend from the basic Chernoff bound is to consider non-independent assumptions for random variables. By Gillman~\cite{gillman1998chernoff} and its refinement works~\cite{chung2012chernoff, rao2017sharp}, they changed the independence assumption to Markov dependence and we summarize their works as follows. We are given $\mathfrak{G}$ as a regular $\lambda$-expander graph with vertex set $\mathfrak{V}$, and $g: \mathfrak{V} \rightarrow \mathbb{C}$ as a bounded function. Suppose $v_1, v_2\cdots, v_{\kappa}$ is a stationary random walk of length $\kappa$ on $\mathfrak{G}$, it is shown that:
\begin{eqnarray}\label{eq:Chernoff bound Markov rvs}\mathrm{Pr} \left( \left\Vert \frac{1}{\kappa} \sum\limits_{j=1}^{\kappa} g(v_i) - \mathbb{E}[g] \right\Vert \geq \vartheta \right) \leq 2 \exp( - \Omega (1 - \lambda) \kappa \vartheta^2).
\end{eqnarray}
The value of $\lambda$ is also the second largest eigenvalue of the transition matrix of the underlying graph $\mathfrak{G}$. The bound given in Eq.~\eqref{eq:Chernoff bound Markov rvs} is named as ``Expander Chernoff Bound''. It is natural to generalize Eq.~\eqref{eq:Chernoff bound Markov rvs} to ``Matrix Expander Chernoff Bound''. Wigderson and Xiao in~\cite{wigderson2008derandomizing} began first attempt to obtain partial results of ``Matrix Expander Chernoff Bound'' and the complete solution is given later by Garg et al.~\cite{garg2018matrix}. 
Let $\mathfrak{G} = (\mathfrak{V}, \mathfrak{E})$ be a regular graph whose transition matrix has second largest eigenvalue as $\lambda$, and let $g: \mathfrak{V} \rightarrow \mathbb{C}^{m \times m}$ be a function satisfy following:
\begin{enumerate}
\item For each $v \mathfrak{V}$, $g(v)$ is a Hermitian matrix with $ \left\Vert g(v) \right\Vert \leq 1$;
\item $\sum\limits_{v \in \mathfrak{V}} g(v) = \mathbf{0}$.
\end{enumerate}
Then, for a stationary random walk $v_1, \cdots, v_\kappa$ with $\epsilon \in (0, 1)$, matrix expander Chernoff bound derived by Garg et al.~\cite{garg2018matrix} is expressed as:
\begin{eqnarray}\label{eq:thm:matrix expander Chernoff bound}
\mathrm{P}\left( \lambda_{g, \max} \left(\frac{1}{\kappa} \sum\limits_{j=1}^{\kappa} g(v_j)\right) \geq \epsilon\right) & \leq &  m \exp(- \Omega (1 - \lambda)\kappa \epsilon^2), \nonumber \\
\mathrm{P}\left( \lambda_{g, \min} \left(\frac{1}{\kappa} \sum\limits_{j=1}^{\kappa} g(v_j)\right) \leq - \epsilon\right) & \leq &  m \exp(- \Omega (1 - \lambda)\kappa \epsilon^2),
\end{eqnarray}
where $\lambda_{g, \max}, (\lambda_{g, \min})$ is the largest (smallest) eigenvalue of the summation of $\kappa$ matrices obtained by the mapping $g$. 

The T-product operation between two three order tensors was invented by Kilmer and her collaborators in~\cite{kilmer2011factorization, kilmer2013third} to generalize the traditional matrix product. T-product operation has been shown as an important linear algebra tool in many domains: multilinear algebra~\cite{li2020continuity, zheng2021t, miao2021t, miao2020generalized}, numerical linear algebra~\cite{zhang2018randomized}, signal processing~\cite{zhang2016exact, semerci2014tensor}, machine learning~\cite{settles2007multiple}, image processing~\cite{khalil2021efficient}, computer vision~\cite{zhang2014novel, martin2013order}, low-rank tensor approximation~\cite{xu2013parallel, zhou2017tensor, qi2021tsingular} etc. However, all these applications assume that systems modelled by T-product tensors are deterministic and such assumption is not true and practical in solving T-product tensors associated issues. In recent years, there are more works begin to study random tensors, see~\cite{chang2020convenient},~\cite{chang2021general},~\cite{vershynin2020concentration} and references therein. Chernoff bound for independent sum of random T-product tensors are considered by the author previous work~\cite{chang2021TProdII}. In this work, we will remove the independent assumption of random T-product tensors by considering T-product tensor expander Chernoff bound. We first transform our majorization technique used in ~\cite{chang2021tensorExpander} for random tensors under Einstein product to random tensors under T-product by deriving norm inequalities of T-product tensors, then we apply these inequalities to build our T-product tensor expander Chernoff bound. Although the author also utilizes norm inequalities of T-product tensors in his recent work~\cite{sychang2021TProdBernstein}, the application of these norm inequalities of T-product tensors is different at this work by considering dependent, instead independent, T-product tensors summation. For self-contained presentation purposes, the portion about majorization techniques and norm inequalities of T-product tensors are also presented in this work. The main result of this paper is summarized by the following theorem. 
\begin{restatable}[T-product Tensor Expander Chernoff Bound]{thm}{TproductExpanderThm}\label{thm:tensor expander}
Let $\mathfrak{G} = (\mathfrak{V}, \mathfrak{E})$ be a regular undirected graph whose transition matrix has second eigenvalue $\lambda$, and let $g: \mathfrak{V} \rightarrow \mathbb{R}^{m \times m \times p}$ be a function. We assume following: 
\begin{enumerate}
\item A nonnegative coefficients polynomial raised by the power $s \geq 1$ as $f: x \rightarrow (a_0 + a_1x +a_2 x^2 + \cdots +a_n x^n)^s$ satisfying $f \left(\exp \left( t   \sum\limits_{j=1}^{\kappa} g(v_j) \right) \right) \succeq \exp \left( t f \left(  \sum\limits_{j=1}^{\kappa} g(v_j) \right) \right) $ almost surely;
\item For each $v \in \mathfrak{V}$, $g(v)$ is a symmetric T-product tensor with $f \left(  \sum\limits_{j=1}^{\kappa} g(v_j) \right)$ as TPD T-product tensor;
\item $\left\Vert g(v) \right\Vert \leq r$;
\item  For $\tau \in [\infty, \infty]$, we have constants $C$ and $\sigma$ such that $ \beta_0(\tau)  \leq \frac{C }{\sigma \sqrt{2 \pi}} \exp \left( \frac{-\tau^2}{2 \sigma^2} \right)$. 
\end{enumerate}
Then, we have 
\begin{eqnarray}\label{eq0:thm:tensor expander}
\mathrm{Pr}\left( \left\Vert f \left(  \sum\limits_{j=1}^{\kappa} g(v_j) \right)  \right\Vert_{(k)} \geq\vartheta \right) \leq   ~~~~~~~~~~~~~~~~~~~~~~~~~~~~~~~~~~~~~~~~~~~~~~~~~~~~~~~~~~~~~~~~~~~~~~  \nonumber \\
  \min\limits_{t > 0 } \left\{ (n+1)^{(s-1)} e^{-\vartheta t} \left[ a_0 k  +C \left(  mp+ \sqrt{\frac{ (mp-k)mp}{k}}  \right)\cdot \right. \right.   \nonumber \\
 \left. \left. \sum\limits_{l=1}^n a_l \exp\left( 8 \kappa \overline{\lambda} + 2  (\kappa +8 \overline{\lambda}) lsr t + 2(\sigma (\kappa +8 \overline{\lambda}) lsr )^2 t^2 \right)  \right] \right\}, ~
\end{eqnarray}
where $\overline{\lambda} = 1 - \lambda$. 
\end{restatable}






The rest of this paper is organized as follows. In Section~\ref{sec:T-product Tensors} , we review T-product tensors basic concepts and introduce a powerful scheme about antisymmetric Kronecker product for T-product tensors. In Section~\ref{sec:Multivariate T-product Tensor Norm Inequalities}, we apply a majorization technique to prove T-product tensor norm inequalities. Our main result about the T-product tensor expander Chernoff bound is provided in Section~\ref{sec:T-product Tensor Expander Chernoff Bound}. Finally, concluding remarks are given by Section~\ref{sec:Conclusions}.

\section{T-product Tensors}\label{sec:T-product Tensors} 

In this section, we will introduce fundamental facts about T-product tensors in Section~\ref{sec:T-product Tensor Fundamental Facts}. Several unitarily invariant norms about a T-product tensor are defined in Section~\ref{sec:Unitarily Invariant T-product Tensor Norms}. A powerful scheme about antisymmetric Kronecker product for T-product tensors will be provided by Section~\ref{sec:Antisymmetric  Kronecker Product for T-product Tensors}.

\subsection{T-product Tensor Fundamental Facts}\label{sec:T-product Tensor Fundamental Facts}

For a third order tensor $\mathcal{C} \in \mathbb{R}^{m \times n \times p}$, we define $\mbox{bcirc}$ operation to the tensor $\mathcal{C}$ as:
\begin{eqnarray}
\mbox{bcirc} (\mathcal{C} ) \define \left[
    \begin{array}{ccccc}
       \mathbf{C}^{(1)}  &  \mathbf{C}^{(p)}  &  \mathbf{C}^{(p-1)}  & \cdots  &  \mathbf{C}^{(2)}   \\
       \mathbf{C}^{(2)}  &  \mathbf{C}^{(1)}  &  \mathbf{C}^{(p)}  & \cdots  &  \mathbf{C}^{(3)}   \\
       \vdots  &  \vdots  & \vdots  & \cdots  & \vdots   \\
       \mathbf{C}^{(p)}  &  \mathbf{C}^{(p-1)}  &  \mathbf{C}^{(p-2)}  & \cdots  &  \mathbf{C}^{(1)}   \\
    \end{array}
\right],
\end{eqnarray}
where $\mathbf{C}^{(1)}, \cdots, \mathbf{C}^{(p)} \in \mathbb{C}^{ m \times n}$ are frontal slices of tensor $\mathcal{C}$. The inverse operation of $\mbox{bcirc}$ is denoted as $\mbox{bcirc}^{-1}$ with relation $\mbox{bcirc}^{-1} ( \mbox{bcirc} ( \mathcal{C} )) \define \mathcal{C}$. Another operation to the tensor $\mathcal{C}$ is unfolding, denoted as $\mbox{unfold} (\mathcal{C} )$, which is defined as: 
\begin{eqnarray}
\mbox{unfold} (\mathcal{C} ) \define \left[
    \begin{array}{c}
       \mathbf{C}^{(1)}  \\
       \mathbf{C}^{(2)}  \\
       \vdots    \\
       \mathbf{C}^{(p)}   \\
    \end{array}
\right].
\end{eqnarray}
The inverse operation of $\mbox{unfold}$ is denoted as $\mbox{fold}$ with relation $\mbox{fold}( \mbox{unfold} ( \mathcal{C} )) \define \mathcal{C}$.

The multiplication between two third order tensors, $\mathcal{C} \mathbb{R}^{m \times n \times p}$ and $\mathcal{D} \mathbb{R}^{n \times l \times p} $, is via \emph{T-product} and this multiplication is defined as:
\begin{eqnarray}\label{eq:T prod def}
\mathcal{C} \star \mathcal{D} \define \mbox{fold}\left( \mbox{bcirc}(\mathcal{C}) \cdot \mbox{unfold}(\mathcal{C})\right),
\end{eqnarray}
where $\cdot$ is the standard matrix multiplication. For given third order tensors, if we apply \emph{T-product} to multiply them, we call them \emph{T-product tensors}. A T-product tensor $\mathcal{C} \in \mathbb{R}^{m \times n \times p}$ will be named as square T-product tensor if $m=n$. 

For a symmetric T-product tensor $\mathcal{C} \in \mathbb{R}^{m \times m \times p}$, we define Hermitian transpose of $\mathcal{C}$, denoted by $\mathcal{C}^{\mathrm{H}}$ , as 
\begin{eqnarray}\label{eq:Hermitian Transpose Def}
\mathcal{C}^{\mathrm{H}} = \mbox{bcirc}^{-1}(     (\mbox{bcirc}(\mathcal{C}))^{\mathrm{H}}  ). 
\end{eqnarray}
And a tensor $\mathcal{D} \in \mathbb{C}^{m \times m \times p }$ is called a Hermitian T-product tensor if $ \mathcal{D}^{\mathrm{H}}  = \mathcal{D}$. Similarly, for a symmetric T-product tensor $\mathcal{C} \in \mathbb{R}^{m \times m \times p}$, we define \emph{transpose} of $\mathcal{C}$, denoted by $\mathcal{C}^{T}$ , as 
\begin{eqnarray}\label{eq:Transpose Def}
\mathcal{C}^{T} = \mbox{bcirc}^{-1}(     (\mbox{bcirc}(\mathcal{C}))^{T}  ). 
\end{eqnarray}
And a tensor $\mathcal{D} \in \mathbb{R}^{m \times m \times p }$ is called a symmetric T-product tensor if $ \mathcal{D}^{T}  = \mathcal{D}$. 

The identity tensor $\mathcal{I}_{m,m,p} \in \mathbb{R}^{m \times m \times p }$ can be defined as:
\begin{eqnarray}\label{eq:I_mmp def}
\mathcal{I}_{m,m,p} = \mbox{bcirc}^{-1}(  \mathbf{I}_{mp}  ),
\end{eqnarray}
where $\mathbf{I}_{mp} $ is the identity matrix in $\mathbb{R}^{mp \times mp}$. For a square T-product tensor, $
\mathcal{C} \in \mathbb{R}^{m \times m \times p }$, we say that $\mathcal{C}$ is nonsingular if it has an inverse tensor $\mathcal{D} \in \mathbb{R}^{m \times m \times p } $ such that 
\begin{eqnarray}\label{eq:inverse tensor def}
\mathcal{C} \star \mathcal{D} = \mathcal{D} \star \mathcal{C} = \mathcal{I}_{m, m, p}.
\end{eqnarray}
A zero tensor, denoted as $\mathcal{O}_{mnp} \in \mathbb{C}^{m \times n \times p}$, is a tensor that all elements inside the tensor as $0$. 

For any circular matrix $\mathbf{C} \in \mathbb{R}^{m \times m}$, it can be diagonalized with the normalized Discrete Fourier Transform (DFT) marix, i.e., $\mathbf{C} = \mathbf{F}^{\mathrm{H}}_m \mathbf{D} \mathbf{F}_m   $, where $\mathbf{F}_m$ is the Fourier matrix of size $m \times m$ defined as 
\begin{eqnarray}\label{eq:DFT def}
\mathbf{F}_m \define  \frac{1}{(m \times p)} \left[
    \begin{array}{ccccc}
       1 &  1 & 1 & \cdots  & 1  \\
       1 &  \omega  &  \omega^2  & \cdots  & \omega^{(m-1)}  \\
       \vdots  &  \vdots  & \vdots  & \cdots  & \vdots   \\
       1 &   \omega^{(m-1)} &  \omega^{2(m-1)} & \cdots  &   \omega^{(m-1)(m-1)}   \\
    \end{array}
\right],
\end{eqnarray}
where $\omega = \exp(\frac{2 \pi \iota}{m})$ with $\iota^2 = -1$. This DFT matrix can also be used to diagonalize a T-product tensor as~\cite{kilmer2013third}
\begin{eqnarray}\label{eq:block diagonalized format}
\mbox{bcirc}(\mathcal{C}) = \left( \mathbf{F}^{\mathrm{H}}_m \otimes \mathbf{I}_m \right) \mbox{Diag}\left( \mathbf{C}_i: i \in \{1, \cdots, m \} \right)  \left( \mathbf{F}_m \otimes \mathbf{I}_m \right),
\end{eqnarray}
where $ \otimes$ is Kronecker Product and $ \mbox{Diag}\left( \mathbf{C}_i: i \in \{1, \cdots, m \} \right) \in \mathbb{C}^{mp \times mp}$ is a diagonal block matrix with the $i$-th diagonal block as the matrix $\mathbf{A}_i$. 

The inner product between two T-product tensors $\mathcal{C} \in \mathbb{C}^{m \times n \times p}$ and $\mathcal{D} \in \mathbb{C}^{m \times n \times p}$ is defined as:
\begin{eqnarray}
\langle \mathcal{C}, \mathcal{D} \rangle = \sum\limits_{i, j, k} c^{\ast}_{i, j, k}d_{i, j, k},
\end{eqnarray}
where $\ast$ is the complex conjugate operation.

We say that a symmetric T-product tensor $\mathcal{C} \in \mathbb{R}^{m \times m \times p}$ is a T-positive definite (TPD) tensor if we have 
\begin{eqnarray}\label{eq: TPD def}
\langle \mathcal{X}, \mathcal{C} \star \mathcal{X} \rangle > 0,
\end{eqnarray}
holds for any non-zero T-product tensor $\mathcal{X} \in \mathbb{R}^{m \times 1 \times p}$. Also, we said that a symmetric T-product tensor is a T-positive semidefinite (TPSD) tensor if we have 
\begin{eqnarray}\label{eq: TPSD def}
\langle \mathcal{X}, \mathcal{C} \mathcal{X} \rangle \geq 0,
\end{eqnarray}
holds for any non-zero T-product tensor $\mathcal{X} \in \mathbb{R}^{m \times 1 \times p}$. Given two T-product tensors $\mathcal{C}, \mathcal{D}$, we use $\mathcal{C} \succ (\succeq) \mathcal{D}$ if $\left(\mathcal{C}  - \mathcal{D} \right)$ is a TPSD (TPD) T-product tensor.

We have the following theorem from Theorem 5 in~\cite{zheng2021t}.
\begin{theorem}\label{thm:tensor and matrix pd relation}
If a T-product tensor $\mathcal{C} \in \mathbb{R}^{m \times m \times p}$ can be diagonalized as 
\begin{eqnarray}\label{eq:block diagonalized format}
\mbox{bcirc}(\mathcal{C}) = \left( \mathbf{F}^{\mathrm{H}}_m \otimes \mathbf{I}_m \right) \mbox{Diag}\left( \mathbf{C}_i: i \in \{1, \cdots, m \} \right)  \left( \mathbf{F}_m \otimes \mathbf{I}_m \right),
\end{eqnarray}
where $\mathbf{F}$ is the DFT matrix defined by Eq.~\eqref{eq:DFT def}; then $\mathcal{C}$ is symmetric, TPD (TPSD) if and only if all matrices $\mathbf{C}_i$ are Hermitian, positive definite (positive semidefinite).
\end{theorem}

Let $\mathcal{C} \in \mathbb{R}^{m \times m \times p}$ can be block diagonalized as Eq.~\eqref{eq:block diagonalized format}. Then, a real number $\lambda$ is said to be a \emph{T-eigenvalue} of $\mathcal{C}$, denoted as $\lambda ( \mathcal{C} )$, if it is an eigenvalue of some $\mathbf{C}_i$ for $i \in \{ 1, \cdots, m \}$. The largest and smallest T-eigenvalue of $\mathcal{C}$ are represented by $\lambda_{\max} (\mathcal{C} ) $ and $\lambda_{\min} (\mathcal{C} ) $, respectively. We use $\lambda_{i, j}$ for the $j$-th largest T-eigenvalue of the matrix $\mathbf{C}_i$.  We also use $\sigma_{i, j}$, named as \emph{T-singular values}, for the $j$-th largest singular values of the matrix $\mathbf{C}_i$. 

We define the T-product tensor \emph{trace} for a tensor $\mathcal{C} =  (c_{ijk}) \in \mathbb{C}^{m \times m \times p}$, denoted by $\mathrm{Tr}(\mathcal{C})$, as following
\begin{eqnarray}\label{eq:trace def}
\mathrm{Tr}(\mathcal{C}) \define \sum\limits_{i=1}^{m}\sum\limits_{k=1}^{p} c_{iik},
\end{eqnarray}
which is the summation of all entries in f-diagonal components. Then, we have the following lemma about trace properties. 
\begin{lemma}\label{lma:T product trace properties}
For any tensors $\mathcal{C}, \mathcal{D} \in \mathbb{C}^{m \times m \times p}$, we have 
\begin{eqnarray}\label{eq:trace linearity}
\mathrm{Tr}(c \mathcal{C} + d \mathcal{D}) = c  \mathrm{Tr}(\mathcal{C}) + d \mathrm{Tr}(\mathcal{D}),
\end{eqnarray}
where $c, d$ are two contants. And, the transpose operation will keep the same trace value, i.e., 
\begin{eqnarray}\label{eq:trace transpose same}
\mathrm{Tr}(\mathcal{C}) = \mathrm{Tr}(\mathcal{C}^{T}).
\end{eqnarray}
Finally, we have 
\begin{eqnarray}\label{eq:trace commutativity}
\mathrm{Tr}(\mathcal{C} \star  \mathcal{D} ) = \mathrm{Tr}(\mathcal{D} \star  \mathcal{C}).
\end{eqnarray}
\end{lemma}
\textbf{Proof:}
Eqs.~\eqref{eq:trace linearity} and~\eqref{eq:trace transpose same} are true from trace definition directly. 

From T-product definition, the $i$-th frontal slice matrix of $\mathcal{D} \star  \mathcal{C}$ is 
\begin{eqnarray}\label{eq:i slice DC}
\mathbf{D}^{(i)} \mathbf{C}^{(1)}  + \mathbf{D}^{(i-1)} \mathbf{C}^{(2)} + \cdots +
\mathbf{D}^{(1)} \mathbf{C}^{(i)} + \mathbf{D}^{(m)} \mathbf{C}^{(i+1) } + \cdots
+ \mathbf{D}^{ (i+1)} \mathbf{C}^{(m)}, 
\end{eqnarray}
similarly, the $i$-th frontal slice matrix of $\mathcal{C} \star  \mathcal{D}$ is 
\begin{eqnarray}\label{eq:i slice CD}
\mathbf{C}^{(i)} \mathbf{D}^{(1)}  + \mathbf{C}^{(i-1)} \mathbf{D}^{(2)} + \cdots +
\mathbf{C}^{(1)} \mathbf{D}^{(i)} + \mathbf{C}^{(m)} \mathbf{D}^{(i+1)} + \cdots
+ \mathbf{C}^{ (i+1)}  \mathbf{D}^{(m)}. 
\end{eqnarray}
Because the matrix trace of Eq.~\eqref{eq:i slice DC} and the matrix trace of Eq.~\eqref{eq:i slice CD} are same for each slice $i$ due to linearity and invariant under cyclic permutations of matrix trace, we have Eq.~\eqref{eq:trace commutativity} by summing over all frontal matrix slices. $\hfill \Box$


Below, we will define the \emph{determinant} of a T-product tensor $\mathcal{C} \in \mathbb{R}^{m \times m \times p}$, represented by $\det ( \mathcal{C} )$, as
\begin{eqnarray}\label{eq:def T prod tensor determinant}
\det ( \mathcal{C} ) &=& \prod\limits_{i=1, j=1}^{i=m, j=p}\lambda_{i, j}. 
\end{eqnarray}

%
%

We have the following theorem from Theorem 6 in~\cite{zheng2021t} about symmetric T-product tensor decomposition. 
\begin{theorem}\label{thm:T eigenvalue decomp}
Every symmetric T-product tensor $\mathcal{C} \in \mathbb{R}^{m \times m \times p}$ can be factored as 
\begin{eqnarray}
\mathcal{C} = \mathcal{U}^{T} \star\mathcal{D} \star \mathcal{U},
\end{eqnarray}
where $\mathcal{U}$ is an orthogonal tensor, i.e., $\mathcal{U}^{T} \star \mathcal{U} = \mathcal{I}_{m, m, p}$, and $\mathcal{D}$ is a F-diagonal tensor, i.e., each frontal slice of $\mathcal{D}$ is a diagonal matrix, such that diagonal entries of $\left( \mathbf{F}_m \otimes \mathbf{I}_m \right) \mbox{bcirc}\left( \mathcal{D} \right) \left( \mathbf{F}^{\mathrm{H}}_m \otimes \mathbf{I}_m \right)$ are T-eigenvalues of $\mathcal{C}$. If $\mathcal{C}$ is a TPD (TPSD) tensor, then all of its T-eigenvalues are positive (nonnegative). 
\end{theorem}

From Theorem~\ref{thm:T eigenvalue decomp} and Lemma~\ref{lma:T product trace properties}, we have the fact that 
\begin{eqnarray}\label{eq:trace is eigen sum}
\mathrm{Tr}(\mathcal{C}) = \sum\limits_i \lambda_i ( \mathcal{C} ).
\end{eqnarray}

%

If a symmetric T-product tensor $\mathcal{C} \in \mathbb{R}^{m \times m \times p}$ can be expressed as the format shown by Eq.~\eqref{eq:block diagonalized format}, the T-eigenvalues of  $\mathcal{C}$ with respect to the matrix $\mathbf{C}_i$ are denoted as $\lambda_{i, k_i}$, where $1 \leq k_i \leq m$, and we assume that 
$\lambda_{i, 1} \geq \lambda_{i, 2} \geq \cdots \geq \lambda_{i, m}$ (including multiplicities). Then, $\lambda_{i, k_i}$ is the $k_i$-th largest T-eigenvalue associated to the matrix $\mathbf{C}_i$. If we sort all T-eigenvalues of $\mathcal{C}$ from the largest one to the smallest one, we use $\tilde{k}$, a smallest integer between 1 to $m \times p$ (inclusive) associated with $p$ given positive integers $k_1, k_2, \cdots, k_p$ that satisfies
\begin{eqnarray}\label{eq1:tilde k}
\lambda_{\tilde{k}} = \min\limits_{1 \leq i \leq m} \lambda_{i, k_i},
\end{eqnarray}
and
\begin{eqnarray}\label{eq2:tilde k}
\lambda_{\tilde{k}}  \geq \lambda_{i, k_i + 1},
\end{eqnarray}
for all $1 \leq i \leq p$. Moreover, we set $\tilde{i}$ from $\lambda_{\tilde{k}} $ as
\begin{eqnarray}\label{eq:tilde i}
\tilde{i}= \arg \min\limits_{i }\left\{  \lambda_{\tilde{k}}= \lambda_{i, k_i} \right\}.
\end{eqnarray}
Then, we will have the following Courant-Fischer theorem for T-product tensors. 


\begin{theorem}\label{thm:Courant-Fischer T-product}
Given a symmetric T-product tensor $\mathcal{C} \in \mathbb{R}^{m \times m \times p}$ and $p$ positive integers $k_1, k_2, \cdots, k_p$ with $1 \leq k_i \leq m$, then we have 
\begin{eqnarray}
\lambda_{\tilde{k}} &=& \max\limits_{\substack{S \in \mathbb{R}^{m \times 1 \times p}\\ \dim(\mathrm{S})  = \{k_1, \cdots, k_p \}   }} \min\limits_{\mathcal{X} \in S } \frac{ \langle \mathcal{X}, \mathcal{C} \star \mathcal{X} \rangle }{ \langle \mathcal{X}, \mathcal{X} \rangle }\nonumber \\ 
 &=&  \min\limits_{\substack{T \in \mathbb{R}^{m \times 1 \times p}\\ \dim(T)  = \{n- k_1, \cdots, n-k_{\tilde{i}-1},  n-k_{\tilde{i}}+1,  n-k_{\tilde{i}+1},\cdots,  n-k_p \}   }} \max\limits_{\mathcal{X} \in T} \frac{ \langle \mathcal{X}, \mathcal{C} \star \mathcal{X} \rangle }{ \langle \mathcal{X}, \mathcal{X} \rangle } 
\end{eqnarray}
where $\lambda_{\tilde{k}}$ and $\tilde{i}$ are defined by Eqs.~\eqref{eq1:tilde k},~\eqref{eq2:tilde k} and~\eqref{eq:tilde i}.
\end{theorem}
\textbf{Proof:}

First, we have to express $ \langle \mathcal{X}, \mathcal{C} \star \mathcal{X} \rangle $ by matrices of $\mathbf{C}_i$ and $\mathbf{X}_i$ through the representation shown by Eq.~\eqref{eq:block diagonalized format}. It is
\begin{eqnarray}\label{eq:thm:Courant-Fischer T-product}
\langle \mathcal{X}, \mathcal{C} \star \mathcal{X} \rangle &=& \frac{1}{p} \langle \mbox{bcirc}(\mathcal{X}), \mbox{bcirc}(\mathcal{C}) \mbox{bcirc}(\mathcal{X})    \rangle \nonumber \\
&=& \frac{1}{p} \mathrm{Tr} \left( \mbox{bcirc}(\mathcal{X})^{\mathrm{H}} \mbox{bcirc}(\mathcal{C}) \mbox{bcirc}(\mathcal{X}) \right)\nonumber \\
&=& \frac{1}{p} \mathrm{Tr} \left(\mathbf{F}^{\mathrm{H}}_p  \mbox{Diag}\left( \mathbf{x}^{\mathrm{H}}_i \mathbf{A}_i  \mathbf{x}_i: i \in \{1,\cdots,p\} \right)\mathbf{F}_p  \right)\nonumber \\
&=& \frac{1}{p} \mathrm{Tr} \left(  \mbox{Diag}\left( \mathbf{x}^{\mathrm{H}}_i \mathbf{A}_i  \mathbf{x}_i: i \in \{1,\cdots,p\} \right)  \right) = \frac{1}{p}\sum\limits_{i=1}^p \mathbf{x}^{\mathrm{H}}_i \mathbf{A}_i  \mathbf{x}_i
\end{eqnarray}

We will just verify the first characterization of $\lambda_{\tilde{k}}$. The other is similar. Let $S_i$ be the projection of $S$ to the space with dimension $k_i$ spanned by $\mathbf{v}_{i, 1}, \cdots, \mathbf{v}_{i, k_i}$, for every $\mathbf{x}_i \in S_{i}$, we can write $\mathbf{x}_i = \sum\limits^{k_i}_{j=1} c_{i, j}  \mathbf{v}_{i, j}$. To show that the value $\lambda_{\tilde{k}}$ is achievable, note that 
\begin{eqnarray}
 \frac{ \langle \mathcal{X}, \mathcal{C} \star \mathcal{X} \rangle }{ \langle \mathcal{X}, \mathcal{X} \rangle }
&=&  \frac{  \frac{1}{p}\sum\limits_{i=1}^p \mathbf{x}^{\mathrm{H}}_i \mathbf{A}_i  \mathbf{x}_i  }{  \frac{1}{p}\sum\limits_{i=1}^p \mathbf{x}^{\mathrm{H}}_i  \mathbf{x}_i  } = \frac{  \sum\limits_{i=1}^p    \sum\limits^{k_i}_{j=1} \lambda_{i, j}    c_{i, j}^{\ast} c_{i, j}          }{   \sum\limits_{i=1}^p    \sum\limits^{k_i}_{j=1}   c_{i, j}^{\ast} c_{i, j}  } \nonumber \\
&\geq & \frac{  \sum\limits_{i=1}^p    \sum\limits^{k_i}_{j=1} \lambda_{\tilde{k}}  c_{i, j}^{\ast} c_{i, j}          }{   \sum\limits_{i=1}^p    \sum\limits^{k_i}_{j=1}   c_{i, j}^{\ast} c_{i, j}  } = \lambda_{\tilde{k}}
\end{eqnarray}
To verify that this is the maximum, let $T_{\tilde{i}}$ be the projection of $T$ to the space with dimension $k_{\tilde{i}}$ with dimension $n- k_{\tilde{i}} + 1$, then the intersection of $S$ and $T_{\tilde{i}}$ is not empty. We have
\begin{eqnarray}
\min\limits_{\mathcal{X} \in S } \frac{ \langle \mathcal{X}, \mathcal{C} \star \mathcal{X} \rangle }{ \langle \mathcal{X}, \mathcal{X} \rangle } &\leq & \min\limits_{\mathcal{X} \in S \cap T} \frac{ \langle \mathcal{X}, \mathcal{C} \star \mathcal{X} \rangle }{ \langle \mathcal{X}, \mathcal{X} \rangle }.
\end{eqnarray}
Any such $\mathbf{x}_{\tilde{i}} \in S \cap T_{\tilde{i}}$ can be expressed as $\mathbf{x}_{\tilde{i}} =  \sum\limits^{m}_{ j=k_{\tilde{i} }} c_{\tilde{i}, j}  \mathbf{v}_{\tilde{i} j}$, and any $i$ for $i \neq \tilde{i}$, we have $\mathbf{x}_{i} \in S \cap T_{i}$ expressed as $\mathbf{x}_{ i } =  \sum\limits^{m}_{ j=k_{i} + 1} c_{i, j}  \mathbf{v}_{i, j}$. Then, we have 
\begin{eqnarray}
 \frac{ \langle \mathcal{X}, \mathcal{C} \star \mathcal{X} \rangle }{ \langle \mathcal{X}, \mathcal{X} \rangle }
&=&  \frac{  \frac{1}{p}\sum\limits_{i=1}^p \mathbf{x}^{\mathrm{H}}_i \mathbf{A}_i  \mathbf{x}_i  }{  \frac{1}{p}\sum\limits_{i=1}^p \mathbf{x}^{\mathrm{H}}_i  \mathbf{x}_i  } = \frac{  \sum\limits_{i=1}^p     \sum\limits^{m}_{ \substack{j=k_i + 1; i \neq \tilde{i} \\ j=k_{\tilde{i}};  i = \tilde{i}    }}   \lambda_{i, j}    c_{i, j}^{\ast} c_{i, j}          }{    \sum\limits_{i=1}^p    \sum\limits^{m}_{ \substack{j=k_i + 1; i \neq \tilde{i} \\ j=k_{\tilde{i}};  i = \tilde{i}    }   } c_{i, j}^{\ast} c_{i, j}   } \nonumber \\
&\leq & \frac{  \sum\limits_{i=1}^p    \sum\limits^{m}_{ \substack{j=k_i + 1; i \neq \tilde{i} \\ j=k_{\tilde{i}};  i = \tilde{i}    }   } \lambda_{\tilde{k}}  c_{i, j}^{\ast} c_{i, j}          }{   \sum\limits_{i=1}^p    \sum\limits^{m}_{ \substack{j=k_i + 1; i \neq \tilde{i} \\ j=k_{\tilde{i}};  i = \tilde{i}    }   } c_{i, j}^{\ast} c_{i, j}  } = \lambda_{\tilde{k}}.
\end{eqnarray}
Therefore, for all subspaces $S$ of dimensions $\{k_1, \cdots, k_p\}$, we have $\min\limits_{\mathcal{X} \in S} \frac{ \langle \mathcal{X}, \mathcal{C} \star \mathcal{X} \rangle }{ \langle \mathcal{X}, \mathcal{X} \rangle } \leq \lambda_{\tilde{k}}$
$\hfill \Box$


\subsection{Unitarily Invariant T-product Tensor Norms}\label{sec:Unitarily Invariant T-product Tensor Norms}


Let us represent the T-eigenvalues of a symmetric T-product tensor $\mathcal{H} \in \mathbb{R}^{m \times m \times p} $ in decreasing order by the vector $\vec{\lambda}(\mathcal{H}) = (\lambda_1(\mathcal{H}), \cdots, \lambda_{m \times p}(\mathcal{H}))$, where $m \times p$ is the total number of T-eigenvalues. We use $\mathbb{R}_{\geq 0} (\mathbb{R}_{> 0})$ to represent a set of nonnegative (positive) real numbers. Let $\left\Vert \cdot \right\Vert_{\rho}$ be a unitarily invariant tensor norm, i.e., $\left\Vert \mathcal{H}\star \mathcal{U}\right\Vert_{\rho} = \left\Vert \mathcal{U}\star \mathcal{H}\right\Vert_{\rho} = \left\Vert \mathcal{H}\right\Vert_{\rho} $,  where $\mathcal{U}$ is any unitary tensor. Let $\rho : \mathbb{R}_{\geq 0}^{m \times p} \rightarrow \mathbb{R}_{\geq 0}$ be the corresponding gauge function that satisfies H$\ddot{o}$lder’s inequality so that 
\begin{eqnarray}\label{eq:def gauge func and general unitarily invariant norm}
\left\Vert \mathcal{H} \right\Vert_{\rho} = \left\Vert |\mathcal{H}| \right\Vert_{\rho} = \rho(\vec{\lambda}( | \mathcal{H} | ) ),
\end{eqnarray}
where $ |\mathcal{H}|  \define \sqrt{\mathcal{H}^H \star \mathcal{H}} $. The bijective correspondence between symmetric gauge functions on $\mathbb{R}_{\geq 0}^{m \times p}$ and unitarily invariant norms is due to von Neumann~\cite{fan1955some}. 

Several popular norms can be treated as special cases of unitarily invariant tensor norm. The first one is Ky Fan like $k$-norm~\cite{fan1955some} for tensors. For $k \in \{1,2,\cdots,m \times p\}$, the Ky Fan $k$-norm~\cite{fan1955some} for tensors  $\mathcal{H}  \mathbb{R}^{m \times m \times p} $, denoted as $\left\Vert \mathcal{H}\right\Vert_{(k)}$, is defined as:
\begin{eqnarray}\label{eq: Ky Fan k norm for tensors}
\left\Vert \mathcal{H}\right\Vert_{(k)} \define \sum\limits_{i=1}^{k} \lambda_i(  |\mathcal{H}|  ).
\end{eqnarray}
If $k=1$,  the Ky Fan $k$-norm for tensors is the tensor operator norm, denoted as $ \left\Vert \mathcal{H} \right\Vert$. The second one is Schatten $p$-norm for tensors, denoted as $\left\Vert \mathcal{H}\right\Vert_{p}$, is defined as:
\begin{eqnarray}\label{eq: Schatten p norm for tensors}
\left\Vert \mathcal{H}\right\Vert_{p} \define (\mathrm{Tr}|\mathcal{H}|^p )^{\frac{1}{p}},
\end{eqnarray}
where $ p \geq 1$. If $p=1$, it is the trace norm. 

Following inequality is the extension of H\"{o}lder inequality to gauge function $\rho$ which will be used later to prove majorization relations. 
\begin{lemma}\label{lma:Holder inquality for gauge function}
For $n$ nonnegative real vectors with the dimension $r$, i.e., $\mathbf{b}_i = (b_{i_1}, \cdots, b_{i_r}) \in \mathbb{R}_{\geq 0}^{r}$, and $\alpha > 0$ with $\sum\limits_{i=1}^n \alpha_i = 1$, we have 
\begin{eqnarray}\label{eq1:lma:Holder inquality for gauge function}
\rho\left( \prod\limits_{i=1}^n b_{i_1}^{\alpha_i},  \prod\limits_{i=1}^n b_{i_2}^{\alpha_i}, \cdots,  \prod\limits_{i=1}^n b_{i_r}^{\alpha_i}  \right) \leq  \prod\limits_{i=1}^n \rho(\mathbf{b}_i)^{\alpha_i} 
\end{eqnarray}
\end{lemma}
\textbf{Proof:}
This proof is based on mathematical induction. The base case for $n=2$ has been shown by Theorem IV.1.6 from~\cite{bhatia2013matrix}. 

We assume that Eq.~\eqref{eq1:lma:Holder inquality for gauge function} is true for $n=m$, where $m > 2$. Let $\odot$ be the component-wise product (Hadamard product) between two vectors.  Then, we have 
\begin{eqnarray}\label{eq2:lma:Holder inquality for gauge function}
\rho\left( \prod\limits_{i=1}^{m+1} b_{i_1}^{\alpha_i},  \prod\limits_{i=1}^{m+1} b_{i_2}^{\alpha_i}, \cdots,  \prod\limits_{i=1}^{m+1} b_{i_r}^{\alpha_i}  \right) = 
\rho\left( \odot_{i=1}^{m+1} \mathbf{b}_i^{\alpha_i}  \right),
\end{eqnarray}
where $\odot_{i=1}^{m+1} \mathbf{b}_i^{\alpha_i}$ is defined as $\left( \prod\limits_{i=1}^{m+1} b_{i_1}^{\alpha_i},  \prod\limits_{i=1}^{m+1} b_{i_2}^{\alpha_i}, \cdots,  \prod\limits_{i=1}^{m+1} b_{i_r}^{\alpha_i}  \right)$ with $\mathbf{b}_i^{\alpha_i} \define (b_{i_1}^{\alpha_i}, \cdots, b_{i_r}^{\alpha_i})$. Under such notations, Eq.~\eqref{eq2:lma:Holder inquality for gauge function} can be bounded as  
\begin{eqnarray}\label{eq3:lma:Holder inquality for gauge function}
\rho\left( \odot_{i=1}^{m+1} \mathbf{b}_i^{\alpha_i}  \right) &= &
\rho\left( \left( \odot_{i=1}^{m} \mathbf{b}_i^{\frac{\alpha_i}{ \sum\limits_{j=1}^m \alpha_j }  } \right)^{\sum\limits_{j=1}^m \alpha_j } \odot \mathbf{b}_{m+1}^{\alpha_{m+1}}\right) \nonumber \\
& \leq & \left[ \rho^{\sum\limits_{j=1}^m \alpha_j } \left( \odot_{i=1}^{m} \mathbf{b}_i^{\frac{\alpha_i}{ \sum\limits_{j=1}^m \alpha_j }  }  \right)  \right] \cdot \rho( \mathbf{b}_{m+1})^{\alpha_{m+1}} \leq  \prod\limits_{i=1}^{m+1} \rho(\mathbf{b}_i)^{\alpha_i}. 
\end{eqnarray}
By mathematical induction, this lemma is proved. $\hfill \Box$

\subsection{Antisymmetric Kronecker Product for T-product Tensors}\label{sec:Antisymmetric  Kronecker Product for T-product Tensors}

In this section, we will discuss a machinery of antisymmetric Kronecker product for T-product tensors and this scheme will be used later for log-majorization results. Let $\mathfrak{H}$ be an $m \times p$-dimensional Hilbert space. For each $k \in \mathbb{N}$, let $\mathfrak{H}^{\otimes k}$ denote the $k$-fold Kronecker product of $\mathfrak{H}$, which is the $(m \times p)^k$-dimensional Hilbert space with respect to the inner product defined by
\begin{eqnarray}
\langle \mathbf{X}_1 \otimes \cdots \otimes \mathbf{X}_k,  \mathbf{Y}_1 \otimes \cdots \otimes \mathbf{Y}_k \rangle \define \prod\limits_{i=1}^k \langle \mathbf{X}_i, \mathbf{Y}_i \rangle.
\end{eqnarray}
For $\mathbf{X}_1, \cdots, \mathbf{X}_k \in \mathfrak{H}$, we define $\mathbf{X}_1 \wedge \cdots \wedge  \mathbf{X}_k \in \mathfrak{H}^{\otimes k}$ by
\begin{eqnarray}
\mathbf{X}_1 \wedge \cdots \wedge  \mathbf{X}_k \define \frac{1}{\sqrt{k!}} \sum\limits_{\sigma} (\mbox{sgn} \sigma) \mathbf{X}_{\sigma(1)} \otimes \cdots \otimes  \mathbf{X}_{\sigma(k)}, 
\end{eqnarray}
where $\sigma$ runs over all permutations on $\{1, 2, \cdots, k\}$ and $\mbox{sgn} \sigma = \pm 1$ depending on $\sigma$ is even or odd. The subspace of $\mathfrak{H}^{\otimes k}$ spanned by $\{\mathbf{X}_1 \wedge \cdots \wedge  \mathbf{X}_k \}$, where $\mathbf{X}_i \in \mathfrak{H}$, is named as $k$-fold antisymmetric Kronecker product of $\mathfrak{H}$ and represented by $\mathfrak{H}^{\wedge k}$. 

For each $\mathcal{C} \in \mathbb{R}^{m \times m \times p}$ and $k \in \mathbb{N}$, the $k$-fold Kronecker product $\mathcal{C}^{\otimes k} \in \mathbb{R}^{m^k \times m^k \times p^k} $ is given by
\begin{eqnarray}
\mathcal{C}^{\otimes k} \star \left( \mathbf{X}_1 \otimes \cdots \otimes \mathbf{X}_k \right) \define 
\left( \mathcal{C} \star \mathbf{X}_1 \right) \otimes  \cdots \otimes \left( \mathcal{C} \star \mathbf{X}_k \right).
\end{eqnarray}
Because $\mathfrak{H}^{\wedge k}$ is invariant for $\mathcal{C}^{\otimes k}$, the antisymmetric Kronecker product of $\mathcal{C}^{\wedge k}$ of $\mathcal{C}$ can be defined as $\mathcal{C}^{\wedge k} = \mathcal{C}^{\otimes}|_{\mathfrak{H}^{\wedge k}}$, then we have
\begin{eqnarray}
\mathcal{C}^{\wedge k} \star \left( \mathbf{X}_1 \wedge \cdots \wedge \mathbf{X}_k \right) 
= \left( \mathcal{C} \star \mathbf{X}_1 \right)  \wedge  \cdots \wedge  \left( \mathcal{C} \star \mathbf{X}_k \right).
\end{eqnarray}

We will provide the following lemmas about antisymmetric Kronecker product.
\begin{lemma}\label{lma:antisymmetric tensor product properties}
Let $\mathcal{A}, \mathcal{B}, \mathcal{C},  \mathcal{E} \in \mathbb{R}^{m \times m \times p}$ be T-product tensors , for any $k \in \{1,2,\cdots,m \times p\}$, we have 
\begin{enumerate}
	\item $(\mathcal{A}^{\wedge k})^\mathrm{T} = (\mathcal{A}^\mathrm{T} )^{\wedge k}$.
	\item $(\mathcal{A}^{\wedge k}) \star (\mathcal{B}^{\wedge k})= (\mathcal{A}\star \mathcal{B})^{\wedge k}$. 
	\item If $\lim\limits_{i \rightarrow \infty} \left\Vert \mathcal{A}_i -  \mathcal{A} \right\Vert \rightarrow 0$ , then $\lim\limits_{i \rightarrow \infty} \left\Vert \mathcal{A}^{\wedge k}_i -  \mathcal{A}^{\wedge k} \right\Vert \rightarrow 0$.
	\item If $\mathcal{C} \succeq \mathcal{O}$ (zero tensor), then $\mathcal{C}^{\wedge k} \succeq \mathcal{O}$ and $(\mathcal{C}^p)^{\wedge k} = (\mathcal{C}^{\wedge k})^p$ for all $p \in \mathbb{R}_{> 0 }$.
    \item  $|\mathcal{A}|^{\wedge k} = | \mathcal{A}^{\wedge k}|$.
    \item  If $\mathcal{E} \succeq \mathcal{O}$ and $\mathcal{E}$ is invertible,  $(\mathcal{E}^z)^{\wedge k} = (\mathcal{E}^{\wedge k})^z$ for all $z \in \mathbb{E}$.
    \item  $\left\Vert \mathcal{E}^{\wedge k} \right\Vert = \prod\limits_{i=1}^{k} \lambda_i ( | \mathcal{E} |)$.
\end{enumerate}
\end{lemma}
\textbf{Proof:}
Items 1 and 2 are the restrictions of the associated relations $(\mathcal{A}^H)^{\otimes k} = (\mathcal{A}^{\otimes k})^H$ and $(\mathcal{A} \star \mathcal{B})^{\otimes k} = (\mathcal{A}^{\otimes k})\star (\mathcal{B}^{\otimes k})$ to $\mathfrak{H}^{\wedge k}$. The item 3 is true since, if $\lim\limits_{i \rightarrow \infty} \left\Vert \mathcal{A}_i -  \mathcal{A} \right\Vert \rightarrow 0$, we have $\lim\limits_{i \rightarrow \infty} \left\Vert \mathcal{A}^{\otimes k}_i -  \mathcal{A}^{\otimes k} \right\Vert \rightarrow 0$ and the asscoaited restrictions of $\mathcal{A}_i^{\otimes k}, \mathcal{A}^{\otimes k}$ to the antisymmetric subspace $\mathfrak{H}^k$. 

For the item 4, if $\mathcal{C} \succeq \mathcal{O}$, then we have $\mathcal{C}^{\wedge k} = ((\mathcal{C}^{1/2})^{\wedge k})^H \star  ((\mathcal{C}^{1/2})^{\wedge k}) \succeq   \mathcal{O}$ from items 1 and 2. If $p$ is ratonal, we have  $(\mathcal{C}^p)^{\wedge k} = (\mathcal{C}^{\wedge k})^p$  from the item 2, and the equality $(\mathcal{C}^p)^{\wedge k} = (\mathcal{C}^{\wedge k})^p$ is also true for any $p > 0$ if we apply the item 3 to approximate any irrelational numbers by rational numbers.

Because we have 
\begin{eqnarray}
|\mathcal{A}|^{\wedge k} =  \left( \sqrt{\mathcal{A}^H \mathcal{A}}\right)^{\wedge k}  =   \sqrt{ (\mathcal{A}^{\wedge k})^H \mathcal{A}^{\wedge k}  }=  | \mathcal{A}^{\wedge k}|,
\end{eqnarray}
from items 1, 2 and 4, so the item 5 is valid. 

For item 6, if $z  < 0$, item 6 is true for all $z \in \mathbb{R}$ by applying the item 4 to $\mathcal{E}^{-1}$. Since we can apply the definition $\mathcal{E}^{z} \define \exp(z \ln \mathcal{E})$ to have
\begin{eqnarray}
\mathcal{C}^p &=& \mathcal{E}^z~~\leftrightarrow~~\mathcal{C} = \exp\left(\frac{z}{p} \ln \mathcal{E} \right),
\end{eqnarray}
where $\mathcal{C} \succeq \mathcal{O}$. The general case of any $z \in \mathbb{C}$ is also true by applying the item 4 to $\mathcal{C} = \exp(\frac{z}{p} \ln \mathcal{E})$. 

For the item 7 proof, it is enough to prove the case that $\mathcal{E} \succeq \mathcal{O}$ due to the item 5. Then, from Theorem~\ref{thm:T eigenvalue decomp}, there exists a set of orthogonal tensors $\{\mathcal{U}_1, \cdots, \mathcal{U}_r\}$ such that $ | \mathcal{E} | \star \mathcal{U}_i = \lambda_i   \mathcal{U}_i$ for $1 \leq i \leq m \times p$. We then have 
\begin{eqnarray}
 | \mathcal{E} |^{\wedge k} \left(\mathcal{U}_{i_1} \wedge \cdots \wedge \mathcal{U}_{i_k}\right)
&=&  | \mathcal{E} |\star  \mathcal{U}_{i_1} \wedge \cdots \wedge  | \mathcal{E} |\star  \mathcal{U}_{i_k}  \nonumber \\
&=& \left( \prod\limits_{i=1}^{k} \lambda_i (  | \mathcal{E} |)  \right) \mathcal{U}_{i_1} \wedge \cdots \wedge \mathcal{U}_{i_k},
\end{eqnarray}
where $1 \leq i_1 < i_2 < \cdots < i_k \leq m \times p$. Hence, $\left\Vert  | \mathcal{E} |^{\wedge k} \right\Vert = \prod\limits_{i=1}^{k} \lambda_i (  | \mathcal{E} |)$. 
$\hfill \Box$

\section{Multivariate T-product Tensor Norm Inequalities}\label{sec:Multivariate T-product Tensor Norm Inequalities}

In this section, we will begin with the introduction of majorization techniques in Section~\ref{sec:Majorization Basis}.
Then, the majorization with integral average and log-majorization with integral average will be introduced by Section~\ref{sec:Majorization wtih Integral Average} and Section~\ref{sec:Log-Majorization wtih Integral Average}. These majorization results will be used to prove T-product tensor norm inequalities in Section~\ref{sec:T-product Tensor Norm Inequalities by Majorization}.


\subsection{Majorization Basis}\label{sec:Majorization Basis} 

In this subsection, we will discuss majorization and several lemmas about majorization which will be used at later proofs. 

Let $\mathbf{x} = [x_1, \cdots,x_r] \in \mathbb{R}^{m \times p}, \mathbf{y} = [y_1, \cdots,y_r] \in \mathbb{R}^{m \times p}$ be two vectors with following orders among entries $x_1 \geq \cdots \geq x_r$ and $y_1 \geq \cdots \geq y_r$, \emph{weak majorization} between vectors $\mathbf{x}, \mathbf{y}$, represented by $\mathbf{x} \prec_{w} \mathbf{y}$, requires following relation for  vectors $\mathbf{x}, \mathbf{y}$:
\begin{eqnarray}\label{eq:weak majorization def}
\sum\limits_{i=1}^k x_i \leq \sum\limits_{i=1}^k y_i,
\end{eqnarray}
where $k \in \{1,2,\cdots,r\}$. \emph{Majorization} between vectors $\mathbf{x}, \mathbf{y}$, indicated by $\mathbf{x} \prec \mathbf{y}$, requires following relation for vectors $\mathbf{x}, \mathbf{y}$:
\begin{eqnarray}\label{eq:majorization def}
\sum\limits_{i=1}^k x_i &\leq& \sum\limits_{i=1}^k y_i,~~\mbox{for $1 \leq k < r$;} \nonumber \\
\sum\limits_{i=1}^{m \times p} x_i &=& \sum\limits_{i=1}^{m \times p} y_i,~~\mbox{for $k = r$.}
\end{eqnarray}

For $\mathbf{x}, \mathbf{y} \in \mathbb{R}^{m \times p}_{\geq 0}$ such that  $x_1 \geq \cdots \geq x_r$ and $y_1 \geq \cdots \geq y_r$,  \emph{weak log majorization} between vectors $\mathbf{x}, \mathbf{y}$, represented by $\mathbf{x} \prec_{w \log} \mathbf{y}$, requires following relation for vectors $\mathbf{x}, \mathbf{y}$:
\begin{eqnarray}\label{eq:weak log majorization def}
\prod\limits_{i=1}^k x_i \leq \prod\limits_{i=1}^k y_i,
\end{eqnarray}
where $k \in \{1,2,\cdots,r\}$, and \emph{log majorization} between vectors $\mathbf{x}, \mathbf{y}$, represented by $\mathbf{x} \prec_{\log} \mathbf{y}$, requires equality for $k=r$ in Eq.~\eqref{eq:weak log majorization def}. If $f$ is a single variable function, $f(\mathbf{x})$ represents a vector of $[f(x_1),\cdots,f(x_r)]$. From Lemma 1 in~\cite{hiai2017generalized}, we have 
\begin{lemma}\label{lma:Lemma 1 Gen Log Hiai}
(1) For any convex function $f: [0, \infty) \rightarrow [0, \infty)$, if we have $\mathbf{x} \prec \mathbf{y}$, then $f(\mathbf{x}) \prec_{w} f(\mathbf{y})$. \\
(2) For any convex function and non-decreasing $f: [0, \infty) \rightarrow [0, \infty)$, if we have $\mathbf{x} \prec_{w} \mathbf{y}$, then $f(\mathbf{x}) \prec_{w} f(\mathbf{y})$. \\
\end{lemma}

Another lemma is from Lemma 12 in~\cite{hiai2017generalized}, we have 
\begin{lemma}\label{lma:Lemma 12 Gen Log Hiai}
Let $\mathbf{x}, \mathbf{y} \in \mathbb{R}^{m \times p}_{\geq 0}$ such that  $x_1 \geq \cdots \geq x_r$ and $y_1 \geq \cdots \geq y_r$ with $\mathbf{x}\prec_{\log} \mathbf{y}$. Also let $\mathbf{y}_i = [y_{i;1}, \cdots , y_{i;r} ] \in \mathbb{R}^{m \times p}_{\geq 0}$ be a sequence of vectors such that $y_{i;1} \geq \cdots \geq y_{i;r} > 0$ and $\mathbf{y}_i \rightarrow \mathbf{y}$ as $i \rightarrow \infty$. Then, there exists $i_0 \in \mathbb{N}$ and $\mathbf{x}_i  = [x_{i;1}, \cdots , x_{i;r} ] \in \mathbb{R}^{m \times p}_{\geq 0}$ for $i \geq i_0$ such that $x_{i;1} \geq \cdots \geq x_{i;r} > 0$, $\mathbf{x}_i \rightarrow \mathbf{x}$ as $i \rightarrow \infty$, and 
\begin{eqnarray}
\mathbf{x}_i \prec_{\log} \mathbf{y}_i \mbox{~~for $i \geq i_0$.} 
\end{eqnarray}
\end{lemma}

For any function $f$ on $\mathbb{R}_{\geq 0}$, the term $f(\mathbf{x}$ is defined as $f(\mathbf{x}) \define (f(x_1), \cdots, f(x_r))$ with conventions $e^{ - \infty} = 0$ and $\log 0 = - \infty$. 

\subsection{Majorization with Integral Average}\label{sec:Majorization wtih Integral Average}

Let $\Omega$ be a $\sigma$-compact metric space and $\nu$ a probability measure on the Borel $\sigma$-field of $\Omega$. Let $\mathcal{C}, \mathcal{D}_\tau \in \mathbb{R}^{m \times m \times p}$ be symmetric T-product tensors. We further assume that tensors $\mathcal{C}, \mathcal{D}_\tau$ are uniformly bounded in their norm for $\tau \in \Omega$. Let $\tau \in\Omega \rightarrow  \mathcal{D}_\tau$ be a continuous function such that $\sup \{\left\Vert  D_{\tau} \right\Vert: \tau \in \Omega  \} < \infty$. For notational convenience, we define the following relation:
\begin{eqnarray}\label{eq:integral eigen vector rep}
\left[ \int_{\Omega} \lambda_1(\mathcal{D}_\tau) d\nu(\tau), \cdots, \int_{\Omega} \lambda_{m \times p}(\mathcal{D}_\tau) d\nu(\tau) \right] \define \int_{\Omega^{m \times p}} \vec{\lambda}(\mathcal{D}_\tau) d\nu^{m \times p}(\tau).
\end{eqnarray}
If $f$ is a single variable function, the notation $f(\mathcal{C})$ represents a tensor function with respect to the tensor $\mathcal{C}$. 

\begin{theorem}\label{thm:weak int average thm 4}
Let $\Omega, \nu, \mathcal{C}, \mathcal{D}_\tau$ be defined as the beginning part of Section~\ref{sec:Majorization wtih Integral Average}, and $f: \mathbb{R} \rightarrow [0, \infty)$ be a non-decreasing convex function, we have following two equivalent statements:
\begin{eqnarray}\label{eq1:thm:weak int average thm 4}
\vec{\lambda}(\mathcal{C}) \prec_w  \int_{\Omega^{m \times p}} \vec{\lambda}(\mathcal{D}_\tau) d\nu^{m \times p}(\tau) \Longleftrightarrow \left\Vert f(\mathcal{C}) \right\Vert_{\rho} \leq 
\int_{\Omega} \left\Vert f(\mathcal{D}_{\tau}) \right\Vert_{\rho}  d\nu(\tau),
\end{eqnarray}
where $\left\Vert \cdot \right\Vert_{\rho}$ is the unitarily invariant norm defined in Eq.~\eqref{eq:def gauge func and general unitarily invariant norm}. 
\end{theorem}
\textbf{Proof:}
We assume that the left statement of Eq.~\eqref{eq1:thm:weak int average thm 4} is true and the function $f$ is a non-decreasing convex function. From Lemma~\ref{lma:Lemma 1 Gen Log Hiai}, we have 
\begin{eqnarray}\label{eq2:thm:weak int average thm 4}
\vec{\lambda}(f (\mathcal{C})) = f (\vec{\lambda}(\mathcal{C})) \prec_w  f \left(\int_{\Omega^{m \times p}} \vec{\lambda}(\mathcal{D}_\tau) d\nu^{m \times p}(\tau) \right).
\end{eqnarray}
From the convexity of $f$, we also have 
\begin{eqnarray}\label{eq3:thm:weak int average thm 4}
f \left(\int_{\Omega^{m \times p} } \vec{\lambda}(\mathcal{D}_\tau) d\nu^{m \times p}  (\tau) \right) \leq \int_{\Omega^{m \times p} } f(\vec{\lambda}(\mathcal{D}_\tau)) d\nu^{m \times p} (\tau) = \int_{\Omega^{m \times p}} \vec{\lambda} ( f(\mathcal{D}_\tau)) d\nu^{m \times p}(\tau).
\end{eqnarray}
Then, we obtain $\vec{\lambda}(f (\mathcal{C}))  \prec_{w} = \int_{\Omega^{m \times p}} \vec{\lambda} ( f(\mathcal{D}_\tau)) d\nu^{m \times p} (\tau)$. By applying Lemma 4.4.2 in~\cite{hiai2010matrix} to both sides of $\vec{\lambda}(f (\mathcal{C}))  \prec_{w} = \int_{\Omega^{m \times p} } \vec{\lambda} ( f(\mathcal{D}_\tau)) d\nu^{m \times p} (\tau)$ with gauge function $\rho$, we obtain 
\begin{eqnarray}\label{eq4:thm:weak int average thm 4}
\left \Vert f(\mathcal{C}) \right\Vert_{\rho} &\leq &\rho \left( \int_{\Omega^{m \times p} } \vec{\lambda} ( f(\mathcal{D}_\tau)) d\nu^{m \times p}  (\tau)  \right)  \nonumber \\
&\leq & \int_{\Omega} \rho(\vec{\lambda} ( f(\mathcal{D}_\tau))) d\nu(\tau) 
= \int_{\Omega} \left\Vert f(\mathcal{D}_\tau) \right\Vert_{\rho} d\nu(\tau).
\end{eqnarray}
Therefore, the right statement of Eq.~\eqref{eq1:thm:weak int average thm 4} is true from the left statement. 

On the other hand, if the right statement of Eq.~\eqref{eq1:thm:weak int average thm 4} is true, we select a function $f \define \max\{x + c, 0\} $, where $c$ is a positive real constant satisfying $\mathcal{C} + c \mathcal{I} \geq \mathcal{O}$, $\mathcal{D}_{\tau} + c \mathcal{I} \geq \mathcal{O}$ for all $\tau \in \Omega$, and tensors $\mathcal{C} + c \mathcal{I}, \mathcal{D}_{\tau} + c \mathcal{I}$. If the Ky Fan $k$-norm at the right statement of Eq.~\eqref{eq1:thm:weak int average thm 4} is applied, we have 
\begin{eqnarray}\label{eq5:thm:weak int average thm 4}
\sum\limits_{i=1}^k (\lambda_i (\mathcal{C}) + c ) \leq  
\sum\limits_{i=1}^k \int_{\Omega} ( \lambda_i (\mathcal{D}_{\tau}) + c ) d\nu(\tau).
\end{eqnarray}
Hence, $\sum\limits_{i=1}^k \lambda_i (\mathcal{C}) \leq  
\sum\limits_{i=1}^k \int_{\Omega} \lambda_i (\mathcal{D}_{\tau}) d\nu(\tau)$, this is the left statement of Eq.~\eqref{eq1:thm:weak int average thm 4}.
$\hfill \Box$

Next theorem will provide a stronger version of Theorem~\ref{thm:weak int average thm 4} by removing weak majorization conditions. 
\begin{theorem}\label{thm:weak int average thm 5}
Let $\Omega, \nu, \mathcal{C}, \mathcal{D}_\tau$ be defined as the beginning part of Section~\ref{sec:Majorization wtih Integral Average}, and $f: \mathbb{R} \rightarrow [0, \infty)$ be a convex function, we have following two equivalent statements:
\begin{eqnarray}\label{eq1:thm:weak int average thm 5}
\vec{\lambda}(\mathcal{C}) \prec  \int_{\Omega^{m \times p}} \vec{\lambda}(\mathcal{D}_\tau) d\nu^{m \times p}(\tau) \Longleftrightarrow \left\Vert f(\mathcal{C}) \right\Vert_{\rho} \leq 
\int_{\Omega} \left\Vert f(\mathcal{D}_{\tau}) \right\Vert_{\rho}  d\nu(\tau),
\end{eqnarray}
where $\left\Vert \cdot \right\Vert_{\rho}$ is the unitarily invariant norm defined in Eq.~\eqref{eq:def gauge func and general unitarily invariant norm}. 
\end{theorem}
\textbf{Proof:}
We assume that the left statement of Eq.~\eqref{eq1:thm:weak int average thm 5} is true and the function $f$ is a convex function. Again, from Lemma~\ref{lma:Lemma 1 Gen Log Hiai}, we have
\begin{eqnarray}\label{eq2:thm:weak int average thm 5}
\vec{\lambda}(f(\mathcal{A})) = f(\vec{\lambda}(\mathcal{A})) \prec_{w}  f \left( \left(\int_{\Omega^{m \times p}} \vec{\lambda}(\mathcal{D}_\tau) d\nu^{m \times p}(\tau)  \right) \right) \leq \int_{\Omega^{m \times p}} f(\vec{\lambda}(\mathcal{D}_\tau)) d\nu^{m \times p}(\tau),
\end{eqnarray}
then, 
\begin{eqnarray}\label{eq3:thm:weak int average thm 5}
\left\Vert f(\mathcal{A}) \right\Vert_{\rho} &\leq & \rho\left( \int_{\Omega^{m \times p}} f(\vec{\lambda}(\mathcal{D}_\tau)) d\nu^{m \times p}(\tau) \right) \nonumber \\
& \leq & \int_{\Omega}\rho \left( f(\vec{\lambda}(\mathcal{D}_\tau)) \right)d\nu (\tau) = 
\int_{\Omega} \left\Vert f( \mathcal{D}_\tau) \right\Vert_{\rho} d\nu (\tau). 
\end{eqnarray}
This proves the right statement of Eq.~\eqref{eq1:thm:weak int average thm 5}. 

Now, we assume that the right statement of Eq.~\eqref{eq1:thm:weak int average thm 5} is true. From Theorem~\ref{thm:weak int average thm 4}, we already have $\vec{\lambda}(\mathcal{C}) \prec_w  \int_{\Omega^{m \times p}} \vec{\lambda}(\mathcal{D}_\tau) d\nu^{m \times p}(\tau)$.  It is enough to prove $\sum\limits_{i=1}^{m \times p} \lambda_i(\mathcal{C}) \geq \int_{\Omega} \sum\limits_{i=1}^{m \times p} \lambda_i(\mathcal{D}_{\tau}) d \nu(\tau)$. We define a function $f \define \max\{c - x, 0\} $, where $c$ is a positive real constant satisfying $\mathcal{C} \leq c \mathcal{I} $, $\mathcal{D}_{\tau} \leq  c \mathcal{I}$ for all $\tau \in \Omega$ and tensors $c \mathcal{I} - \mathcal{C}, c \mathcal{I} - \mathcal{D}_{\tau}$. If the trace norm is applied, i.e., the sum of the absolute value of all eigenvalues of a symmetric T-product tensor, then the right statement of Eq.~\eqref{eq1:thm:weak int average thm 5} becomes
\begin{eqnarray}\label{eq4:thm:weak int average thm 5}
\sum\limits_{i=1}^{m \times p} \lambda_i \left( c\mathcal{I} - \mathcal{C}\right)  \leq \int_{\Omega} 
\sum\limits_{i=1}^{m \times p} \lambda_i \left( c\mathcal{I} - \mathcal{D}_{\tau}\right) d \nu(\tau).
\end{eqnarray}
The desired inequality  $\sum\limits_{i=1}^{m \times p} \lambda_i(\mathcal{C}) \geq \int_{\Omega} \sum\limits_{i=1}^{m \times p} \lambda_i(\mathcal{D}_{\tau}) d \nu(\tau)$ is established. $\hfill \Box$

\subsection{Log-Majorization with Integral Average}\label{sec:Log-Majorization wtih Integral Average}

The purpose of this section is to consider log-majorization issues for unitarily invariant norms of TPSD T-product tensors. In this section, let $\mathcal{C}, \mathcal{D}_\tau \in \mathbb{R}^{m \times m \times p}$ be TPSD T-product tensors with $m \times p$ nonnegative T-eigenvalues by keeping notations with the same definitions as at the beginning of the Section~\ref{sec:Majorization wtih Integral Average}. For notational convenience, we define the following relation for logarithm vector:
\begin{eqnarray}\label{eq:integral eigen log vector rep}
\left[ \int_{\Omega} \log \lambda_1(\mathcal{D}_\tau) d\nu(\tau), \cdots, \int_{\Omega} \log \lambda_{m \times p}(\mathcal{D}_\tau) d\nu(\tau) \right] \define \int_{\Omega^{m \times p}} \log \vec{\lambda}(\mathcal{D}_\tau) d\nu^{m \times p}(\tau).
\end{eqnarray}

\begin{theorem}\label{thm:weak int average thm 7}
Let $\mathcal{C}, \mathcal{D}_\tau$ be TPSD T-product tensors, $f: (0, \infty) \rightarrow [0,\infty)$ be a continuous function such that the mapping $x \rightarrow \log f(e^x)$ is convex on $\mathbb{R}$, and $g: (0, \infty) \rightarrow [0,\infty)$ be a continuous function such that the mapping $x \rightarrow g(e^x)$ is convex on $\mathbb{R}$ , then  we have following three equivalent statements:
\begin{eqnarray}\label{eq1:thm:weak int average thm 7}
\vec{\lambda}(\mathcal{C}) &\prec_{w \log}& \exp  \int_{\Omega^{m \times p}} \log \vec{\lambda}(\mathcal{D}_\tau) d\nu^{m \times p}(\tau);
\end{eqnarray}
\begin{eqnarray}\label{eq2:thm:weak int average thm 7}
\left\Vert f(\mathcal{C}) \right\Vert_{\rho} &\leq &
\exp \int_{\Omega} \log \left\Vert f(\mathcal{D}_{\tau}) \right\Vert_{\rho}  d\nu(\tau);
\end{eqnarray}
\begin{eqnarray}\label{eq3:thm:weak int average thm 7}
\left\Vert g(\mathcal{C}) \right\Vert_{\rho} &\leq &
\int_{\Omega} \left\Vert g(\mathcal{D}_{\tau}) \right\Vert_{\rho}  d\nu(\tau).
\end{eqnarray}
\end{theorem}
\textbf{Proof:}
The roadmap of this proof is to prove equivalent statements between Eq.~\eqref{eq1:thm:weak int average thm 7} and Eq.~\eqref{eq2:thm:weak int average thm 7} first, followed by equivalent statements between Eq.~\eqref{eq1:thm:weak int average thm 7} and Eq.~\eqref{eq3:thm:weak int average thm 7}. 

\textbf{Eq.~\eqref{eq1:thm:weak int average thm 7} $\Longrightarrow$ Eq.~\eqref{eq2:thm:weak int average thm 7}}

There are two cases to be discussed in this part of proof: $\mathcal{C}, \mathcal{D}_\tau$ are TPD tensors, and $\mathcal{C}, \mathcal{D}_\tau$ are TPSD T-product tensors. At the beginning, we consider the case that $\mathcal{C}, \mathcal{D}_\tau$ are TPD tensors.

Since $\mathcal{D}_\tau$ are positive, we can find $\varepsilon > 0$ such that $\mathcal{D}_{\tau} \geq \varepsilon \mathcal{I}$ for all $\tau \in \Omega$. From Eq.~\eqref{eq1:thm:weak int average thm 7}, the convexity of $\log f (e^x)$ and Lemma~\ref{lma:Lemma 1 Gen Log Hiai}, we have 
\begin{eqnarray}
\vec{\lambda} \left( f ( \mathcal{C})\right)  = f \left(\exp \left( \log \vec{\lambda} (\mathcal{C}) \right) \right) &\prec _w &  f \left(\exp    \int_{\Omega^{m \times p}} \vec{\lambda}(\mathcal{D}_\tau) d\nu^{m \times p}(\tau)             \right) \nonumber \\
& \leq &  \exp \left( \int_{\Omega^{m \times p}} \log f \left( \vec{\lambda}(\mathcal{D}_\tau) \right)  d\nu^{m \times p}(\tau) \right).
\end{eqnarray}
Then, from Eq.~\eqref{eq:def gauge func and general unitarily invariant norm}, we obtain
\begin{eqnarray}\label{eq4:thm:weak int average thm 7}
\left\Vert f (\mathcal{C}) \right\Vert _{\rho}
& \leq &  \rho\left(\exp \left( \int_{\Omega^{m \times p}} \log  f \left( \vec{\lambda}(\mathcal{D}_\tau) \right) d\nu^{m \times p}(\tau) \right)   \right).
\end{eqnarray}


From the function $f$ properties, we can assume that $f(x) > 0$ for any $x > 0$. Then, we have 
following bounded and continous maps on $\Omega$: $\tau \rightarrow \log f (\lambda_i (\mathcal{D}_{\tau}))  $ for $i \in \{1,2,\cdots, m \times p \}$, and $\tau \rightarrow \log \left\Vert f (\mathcal{D}_{\tau}) \right\Vert_{\rho}$. Because we have $\nu (\Omega) = 1$ and $\sigma$-compactness of $\Omega$, we have $\tau_{k}^{(n)} \in \Omega$ and $\alpha_{k}^{(n)}$ for $k \in \{1,2,\cdots, n\}$ and $n \in \mathbb{N}$ with $\sum\limits_{k=1}^{n} \alpha_{k}^{(n)} = 1$ such that 
\begin{eqnarray}\label{eq:35}
\int_{\Omega} \log f (\lambda_i ( \mathcal{D}_{\tau} ) ) d \nu (\tau) = \lim\limits_{n \rightarrow \infty} \sum\limits_{k=1}^{n} \alpha_{k}^{(n)}  \log f (\lambda_i (\mathcal{D}_{\tau_k^{(n)} }))   , \mbox{for $i \in \{1,2,\cdots  m \times p \}$};
\end{eqnarray}
and 
\begin{eqnarray}\label{eq:36}
\int_{\Omega} \log \left\Vert f (\mathcal{D}_{\tau}) \right\Vert_{\rho} d \nu (\tau) = \lim\limits_{n \rightarrow \infty} \sum\limits_{k=1}^{n} \alpha_{k}^{(n)}  \log \left\Vert f (\mathcal{D}_{\tau_k^{(n)}    }) \right\Vert_{\rho} .
\end{eqnarray}
By taking the exponential at both sides of Eq.~\eqref{eq:35} and apply the gauge function $\rho$, we have
\begin{eqnarray}\label{eq:37}
\rho \left( \exp \int_{\Omega^{m \times p}} \log f (\vec{\lambda} ( \mathcal{D}_{\tau} ) ) d \nu^{m \times p} (\tau)  \right)= \lim\limits_{n \rightarrow \infty} \rho\left(  \prod\limits_{k=1}^{n}  f  \left( \vec{\lambda} \left(\mathcal{D}_{\tau_k^{(n)} } \right)  \right)^{ \alpha_{k}^{(n)} }  \right).
\end{eqnarray}
Similarly, by taking the exponential at both sides of Eq.~\eqref{eq:36}, we have
\begin{eqnarray}\label{eq:38}
\exp \left( \int_{\Omega} \log \left\Vert f (\mathcal{D}_{\tau}) \right\Vert_{\rho} d \nu (\tau) \right) = \lim\limits_{n \rightarrow \infty} \prod \limits_{k=1}^{n} \left\Vert f \left( \mathcal{D}_{\tau_k^{(n)}    } \right) \right\Vert^{\alpha_{k}^{(n)}}_{\rho} .
\end{eqnarray}
From Lemma~\ref{lma:Holder inquality for gauge function}, we have 
\begin{eqnarray}\label{eq:41}
\rho \left(  \prod\limits_{k=1}^{n}  f  \left( \vec{\lambda} \left(\mathcal{D}_{\tau_k^{(n)} } \right)  \right)^{ \alpha_{k}^{(n)} }  \right) & \leq & \prod \limits_{k=1}^{n} \rho^{  \alpha_{k}^{(n)}    } \left( f \left( \vec{\lambda} \left( \mathcal{D}_{ \tau_{k}^{(n)}  }\right) 
 \right) \right) \nonumber \\
&=& \prod \limits_{k=1}^{n} \rho^{  \alpha_{k}^{(n)}    } \left( \vec{\lambda}  \left( f  \left( \mathcal{D}_{ \tau_{k}^{(n)}  }\right) 
 \right) \right) \nonumber \\
&=& \prod \limits_{k=1}^{n} \left\Vert f  \left( \mathcal{D}_{ \tau_{k}^{(n)}  } \right) \right\Vert_{\rho}^{\alpha_{k}^{(n)}}
\end{eqnarray}

From Eqs.~\eqref{eq:37},~\eqref{eq:38} and~\eqref{eq:41}, we have 
\begin{eqnarray}\label{eq:42}
\rho \left( \exp \int_{\Omega^{m \times p}} \log f (\vec{\lambda} ( \mathcal{D}_{\tau} ) ) d \nu^{m \times p} (\tau)  \right) \leq \exp \int_{\Omega} \log \left\Vert f(\mathcal{D}_{\tau})\right\Vert_{\rho} d \nu(\tau).
\end{eqnarray}
Then, Eq.~\eqref{eq2:thm:weak int average thm 7} is proved from Eqs.~\eqref{eq4:thm:weak int average thm 7} and~\eqref{eq:42}.

Next, we consider that $\mathcal{C}, \mathcal{D}_\tau$ are TPSD T-product tensors. For any $\delta > 0$, we have following log-majorization relation:
\begin{eqnarray}
\prod\limits_{i=1}^k \left( \lambda_i (\mathcal{C}) + \epsilon_{\delta} \right) 
&\leq& \prod\limits_{i=1}^k \exp  \int_{\Omega} \log \left( \lambda_i(\mathcal{D}_\tau) + \delta\right) d \nu (\tau),
\end{eqnarray}
where $\epsilon_{\delta} > 0$ and $k \in \{1,2,\cdots r \}$. Then, we can apply the previous case result about TPD tensors to TPD tensors $\mathcal{C} + \epsilon_{\delta} \mathcal{I}$ and $\mathcal{D}_\tau + \delta \mathcal{I}$, and get 
\begin{eqnarray}\label{eq:46}
\left\Vert f (\mathcal{C}) + \epsilon_{\delta}  \mathcal{I} \right\Vert_{\rho} 
&\leq& \exp \int_{\Omega} \log \left\Vert f (\mathcal{D}_{\tau}) + \delta  \mathcal{I} \right\Vert_{\rho} 
d \nu (\tau)
\end{eqnarray}
As $\delta \rightarrow 0$, Eq.~\eqref{eq:46} will give us Eq.~\eqref{eq2:thm:weak int average thm 7} for TPSD T-product tensors.  

\textbf{Eq.~\eqref{eq1:thm:weak int average thm 7} $\Longleftarrow$ Eq.~\eqref{eq2:thm:weak int average thm 7}}

We consider TPD tensors at first phase by assuming that $\mathcal{D}_{\tau}$ are 
TPD T-product tensors for all $\tau \in \Omega$. We may also assume that the tensor $\mathcal{C}$ is a TPD T-product tensor. Since if this is a TPSD T-product tensor, i.e., some $\lambda_i = 0$, we always have following inequality valid:
\begin{eqnarray}
\prod\limits_{i=1}^k \lambda_i (\mathcal{C}) \leq \prod\limits_{i=1}^k 
\exp \int_{\Omega} \log \lambda_i (\mathcal{D}_{\tau}) d \nu (\tau)
\end{eqnarray}

If we apply $f(x) = x^p$ for $p > 0$ and $\left\Vert \cdot \right\Vert_{\rho}$ as Ky Fan $k$-norm in Eq.~\eqref{eq2:thm:weak int average thm 7}, we have 
\begin{eqnarray}\label{eq:50}
\log \sum\limits_{i=1}^k \lambda^p_i \left(\mathcal{C}\right) \leq \int_{\Omega} \log \sum\limits_{i=1}^k \lambda_i^p\left( \mathcal{D}_{\tau} \right) d \nu (\tau).
\end{eqnarray}
If we add $\log \frac{1}{k}$ and multiply $\frac{1}{p}$ at both sides of Eq.~\eqref{eq:50}, we have 
\begin{eqnarray}\label{eq:51}
\frac{1}{p}\log \left( \frac{1}{k}\sum\limits_{i=1}^k \lambda^p_i \left(\mathcal{C}\right) \right)\leq \int_{\Omega} \frac{1}{p} \log \left( \frac{1}{k}\sum\limits_{i=1}^k \lambda_i^p\left( \mathcal{D}_{\tau} \right) \right) d \nu (\tau).
\end{eqnarray}
From L'Hopital's Rule, if $p \rightarrow 0$, we have 
\begin{eqnarray}\label{eq:52}
\frac{1}{p}\log \left( \frac{1}{k}\sum\limits_{i=1}^k \lambda^p_i \left(\mathcal{C}\right) \right) \rightarrow \frac{1}{k} \sum\limits_{i=1}^k \log \lambda_i (\mathcal{C}),
\end{eqnarray}
and 
\begin{eqnarray}\label{eq:53}
\frac{1}{p}\log \left( \frac{1}{k}\sum\limits_{i=1}^k \lambda^p_i \left(\mathcal{D}_{\tau}\right) \right) \rightarrow \frac{1}{k} \sum\limits_{i=1}^k \log \lambda_i (\mathcal{D}_{\tau}),
\end{eqnarray}
where $\tau \in \Omega$. Appling Eqs.~\eqref{eq:52} and~\eqref{eq:53} into Eq.~\eqref{eq:51} and taking $p \rightarrow 0$, we have 
\begin{eqnarray}
\sum\limits_{i=1}^k \lambda_i (\mathcal{C}) \leq \int_{\Omega} \sum\limits_{i=1}^k 
 \log \lambda_i (\mathcal{D}_{\tau}) d \nu (\tau).
\end{eqnarray}
Therefore, Eq.~\eqref{eq1:thm:weak int average thm 7} is true for TPD tensors. 

For TPSD T-product tensors $\mathcal{D}_{\tau}$, since Eq.~\eqref{eq2:thm:weak int average thm 7} is valid for $\mathcal{D}_{\tau} + \delta \mathcal{I}$ for any $\delta > 0$, we can apply the previous case result about TPD tensors to $\mathcal{D}_{\tau} + \delta \mathcal{I}$ and obtain
\begin{eqnarray}
\prod\limits_{i=1}^k \lambda_i (\mathcal{C}) \leq \prod\limits_{i=1}^k  \exp 
\int_{\Omega}   \log  \left( \lambda_i (\mathcal{D}_{\tau}) + \delta \right) d \nu (\tau),
\end{eqnarray}
where $k \in \{1,2,\cdots, r\}$. Eq.~\eqref{eq1:thm:weak int average thm 7} is still true for TPSD T-product tensors as $\delta \rightarrow 0$.

\textbf{Eq.~\eqref{eq1:thm:weak int average thm 7} $\Longrightarrow$ Eq.~\eqref{eq3:thm:weak int average thm 7}}

If $\mathcal{C}, \mathcal{D}_\tau$ are TPD tensors, and $\mathcal{D}_{\tau} \geq \delta \mathcal{I}$ for all $\tau \in \Omega$. From Eq.~\eqref{eq1:thm:weak int average thm 7}, we have 
\begin{eqnarray}
\vec{\lambda} (\log \mathcal{C}) = \log \vec{\lambda}(\mathcal{C}) \prec_{w}
\int_{\Omega^{m \times p}} \log \vec{\lambda}(\mathcal{D}_{\tau}) d \nu^{m \times p} (\tau) = 
\int_{\Omega^{m \times p}} \vec{\lambda}( \log \mathcal{D}_{\tau}) d \nu^{m \times p} (\tau).
\end{eqnarray}
If we apply Theorem~\ref{thm:weak int average thm 4} to $\log \mathcal{C}$, $\log \mathcal{D}_{\tau}$ with function $f(x) = g(e^x)$, where $g$ is used in Eq.~\eqref{eq3:thm:weak int average thm 7}, Eq.~\eqref{eq3:thm:weak int average thm 7} is implied. 

If $\mathcal{C}, \mathcal{D}_\tau$ are TPSD T-product tensors and any $\delta > 0$, we can find $\epsilon_{\delta} \in (0, \delta)$ to satisfy following:
\begin{eqnarray}\label{eq:45}
\prod\limits_{i=1}^k\left(\lambda_i(\mathcal{C}) + \epsilon_{\delta}\right) \leq 
\prod\limits_{i=1}^k \exp \int_{\Omega}   \log \left( \lambda_i(\mathcal{D}_{\tau}) + \delta  \right) d \nu (\tau).
\end{eqnarray}
Then, from TPD T-product tensor case, we have 
\begin{eqnarray}\label{eq:45-1}
\left\Vert g( \mathcal{C} + \epsilon_{\delta} \mathcal{I} ) \right\Vert_{\rho}
\leq \int_{\Omega} \left\Vert   g( \mathcal{D}_{\tau} + \delta \mathcal{I} )    \right\Vert_{\rho}
d \nu (\tau).
\end{eqnarray}
Eq.~\eqref{eq3:thm:weak int average thm 7} is obtained by taking $\delta \rightarrow 0$ in Eq.~\eqref{eq:45-1}. 

\textbf{Eq.~\eqref{eq1:thm:weak int average thm 7} $\Longleftarrow$ Eq.~\eqref{eq3:thm:weak int average thm 7}}

For $k \in \{1,2,\cdots, r \}$, if we apply $g(x) = \log (\delta + x )$, where $\delta >0$, and Ky Fan $k$-norm in Eq.~\eqref{eq3:thm:weak int average thm 7}, we have 
\begin{eqnarray}
\sum\limits_{i=1}^k \log \left(\delta + \lambda_i \left(\mathcal{C} \right) \right)
\leq \sum\limits_{i=1}^k \int_{\Omega} \log \left( \delta + \lambda_{i}(\mathcal{D}_{\tau}) \right) d \nu (\tau).
\end{eqnarray}
Then, we have following relation as $\delta \rightarrow 0$:
\begin{eqnarray}
\sum\limits_{i=1}^k \log \lambda_i \left(\mathcal{C} \right) 
\leq \sum\limits_{i=1}^k \int_{\Omega} \log  \lambda_{i}(\mathcal{D}_{\tau}) d \nu (\tau).
\end{eqnarray}
Therefore, Eq.~\eqref{eq1:thm:weak int average thm 7} ccan be derived from Eq.~\eqref{eq3:thm:weak int average thm 7}. $\hfill \Box$

Next theorem will extend Theorem~\ref{thm:weak int average thm 7} to non-weak version.

\begin{theorem}\label{thm:int log average thm 10}
Let $\mathcal{C}, \mathcal{D}_\tau$ be TPSD T-product tensors with $\int_{\Omega} \left\Vert \mathcal{D}_{\tau}^{-p}\right\Vert_\rho d \nu (\tau) < \infty$ for any $p > 0$, $f: (0, \infty) \rightarrow [0,\infty)$ be a continuous function such that the mapping $x \rightarrow \log f(e^x)$ is convex on $\mathbb{R}$, and $g: (0, \infty) \rightarrow [0,\infty)$ be a continuous function such that the mapping $x \rightarrow g(e^x)$ is convex on $\mathbb{R}$ , then  we have following three equivalent statements:
\begin{eqnarray}\label{eq1:thm:int average thm 10}
\vec{\lambda}(\mathcal{C}) &\prec_{\log}& \exp  \int_{\Omega^{m \times p}} \log \vec{\lambda}(\mathcal{D}_\tau) d\nu^{m \times p}(\tau);
\end{eqnarray}
\begin{eqnarray}\label{eq2:thm:int average thm 10}
\left\Vert f(\mathcal{C}) \right\Vert_{\rho} &\leq &
\exp \int_{\Omega} \log \left\Vert f(\mathcal{D}_{\tau}) \right\Vert_{\rho}  d\nu(\tau);
\end{eqnarray}
\begin{eqnarray}\label{eq3:thm:int average thm 10}
\left\Vert g(\mathcal{C}) \right\Vert_{\rho} &\leq &
\int_{\Omega} \left\Vert g(\mathcal{D}_{\tau}) \right\Vert_{\rho}  d\nu(\tau).
\end{eqnarray}
\end{theorem}
\textbf{Proof:}

The proof plan is similar to the proof in Theorem~\ref{thm:weak int average thm 7}. We prove the equivalence between Eq.~\eqref{eq1:thm:int average thm 10} and Eq.~\eqref{eq2:thm:int average thm 10} first, then prove the equivalence between Eq.~\eqref{eq1:thm:int average thm 10} and Eq.~\eqref{eq3:thm:int average thm 10}. 

\textbf{Eq.~\eqref{eq1:thm:int average thm 10} $\Longrightarrow$ Eq.~\eqref{eq2:thm:int average thm 10}}

First, we assume that $\mathcal{C}, \mathcal{D}_\tau$ are TPD T-product tensors with $\mathcal{D}_\tau \geq \delta \mathcal{I}$ for all $\tau \in \Omega$. The corresponding part of the proof in Theorem~\ref{thm:weak int average thm 7} about TPD tensors $\mathcal{C}, \mathcal{D}_\tau$ can be applied here. 

For case that $\mathcal{C}, \mathcal{D}_\tau$ are TPSD T-product tensors, we have 
\begin{eqnarray}
\prod\limits_{i=1}^k \lambda_i (\mathcal{C}) \leq \prod\limits_{i=1}^k 
\exp \int_{\Omega} \log \left( \lambda_i (\mathcal{D}_{\tau}) + \delta_n \right)d \nu (\tau), 
\end{eqnarray}
where $\delta_n > 0$ and $\delta_n \rightarrow 0$. Because $\int_{\Omega^{m \times p}} \log \left( \vec{\lambda} (\mathcal{D}_{\tau}) + \delta_n \right) d \nu^{m \times p} (\tau) \rightarrow \int_{\Omega^{m \times p}} \log \vec{\lambda} (\mathcal{D}_\tau)  d \nu^{m \times p} (\tau) $ as $n \rightarrow \infty$, from Lemma~\ref{lma:Lemma 12 Gen Log Hiai}, we can find $\mathbf{a}^{(n)}$ with $n \geq n_0$ such that $a^{(n)}_1 \geq \cdots \geq a^{(n)}_r > 0$, $\mathbf{a}^{(n)} \rightarrow \vec{\lambda}(\mathcal{C})$ and $ \mathbf{a}^{(n)} \prec_{\log}  \exp \int_{\Omega^{m \times p}} \log \vec{\lambda} \left(\mathcal{D}_{\tau} + \delta_n \mathcal{I} \right) d \nu^{m \times p} (\tau)$

Selecting $\mathcal{C}^{(n)}$ with $\vec{\lambda} ( \mathcal{C}^{(n)})  = \mathbf{a}^{(n)} $ and applying TPD tensors case to $\mathcal{C}^{(n)}$ and $\mathcal{D}_{\tau} + \delta_n \mathcal{I}$, we obtain
\begin{eqnarray}\label{eq:74}
\left\Vert f (\mathcal{C}^{(n)}) \right\Vert_{\rho} \leq \exp \int_{\Omega} \log \left\Vert f (\mathcal{D}_{\tau} + \delta_n \mathcal{I}) \right\Vert_{\rho} d \nu (\tau)
\end{eqnarray}
where $n \geq n_0$.

There are two situations for the function $f$ near $0$: $f(0^{+}) < \infty$ and $f(0^{+}) = \infty$. For the case with $f(0^{+}) < \infty$, we have 
\begin{eqnarray}\label{eq:75-1}
\left\Vert f (\mathcal{C}^{(n)} )\right\Vert_{\rho} = \rho( f (\mathbf{a}^{(n)}))
\rightarrow \rho (f ( \vec{\lambda}(\mathcal{C}))) = \left\Vert f (\mathcal{C})\right\Vert_{\rho}, 
\end{eqnarray}
and
\begin{eqnarray}\label{eq:75-2}
\left\Vert f (\mathcal{D}_{\tau} + \delta_n \mathcal{I} )\right\Vert_{\rho} 
\rightarrow \left\Vert f (\mathcal{D}_{\tau})\right\Vert_{\rho}, 
\end{eqnarray}
where $\tau \in \Omega$ and $n \rightarrow \infty$. From Fatou–Lebesgue theorem, we then have 
\begin{eqnarray}\label{eq:76}
\limsup\limits_{n \rightarrow \infty} \int_{\Omega} \log \left\Vert f (\mathcal{D}_{\tau} + \delta_n \mathcal{I} )\right\Vert_{\rho} d \nu (\tau) \leq \int_{\Omega} \log \left\Vert f(\mathcal{D}_{\tau}) \right\Vert_{\rho}.
\end{eqnarray}
By taking $n \rightarrow \infty$ in Eq.~\eqref{eq:74} and using Eqs.~\eqref{eq:75-1},~\eqref{eq:75-2},~\eqref{eq:76}, we have Eq.~\eqref{eq2:thm:int average thm 10} for case that $f(0^{+}) < \infty$.

For the case with $f(0^{+}) = \infty$, we assume that $\int_{\Omega} \log \left\Vert f (\mathcal{D}_{\tau}) \right\Vert_{\rho} d \nu (\tau) < \infty$ (since the inequality in Eq.~\eqref{eq2:thm:int average thm 10} is always true for $\int_{\Omega} \log \left\Vert f (\mathcal{D}_{\tau}) \right\Vert_{\rho} d \nu (\tau) = \infty$). Since $f$ is decreasing on $(0, \epsilon)$ for some $\epsilon > 0$. We claim that the following relation is valid: there are two constants $a, b > 0$ such that 
\begin{eqnarray}\label{eq:77}
a \leq \left\Vert f (\mathcal{D}_{\tau} + \delta_n \mathcal{I}) \right\Vert_{\rho}
\leq \left\Vert f (\mathcal{D}_{\tau}) \right\Vert_{\rho} + b,
\end{eqnarray}
for all $\tau \in \Omega$ and $n \geq n_0$. If Eq.~\eqref{eq:77} is valid and $\int_{\Omega} \log \left\Vert f (\mathcal{D}_{\tau}) \right\Vert_{\rho} d \nu (\tau) < \infty$, from Lebesgue's dominated convergence theorem, we also have Eq.~\eqref{eq2:thm:int average thm 10} for case that $f(0^{+}) = \infty$ by taking $n \rightarrow \infty$ in Eq.~\eqref{eq:74}. 

Below, we will prove the claim stated by Eq.~\eqref{eq:77}. By the uniform boundedness of tensors $\mathcal{D}_{\tau}$, there is a constant $\kappa >0$ such that 
\begin{eqnarray}
0 < \mathcal{D}_{\tau} + \delta_n \mathcal{I} \leq \kappa \mathcal{I},
\end{eqnarray}
where $\tau \in \Omega$ and $ n \geq n_0$. We may assume that $\mathcal{D}_\tau$ is TPD tensors because $\left\Vert f (\mathcal{D}_{\tau}) \right\Vert_{\rho} = \infty$, i.e., Eq.~\eqref{eq:77} being true automatically, when $\mathcal{D}_\tau$ is TPSD T-product tensors. From Theorem~\ref{thm:T eigenvalue decomp}, we have 
\begin{eqnarray}
f(\mathcal{D}_{\tau} + \delta_n \mathcal{I}) &=&  \sum\limits_{i', \mbox{s.t. $\lambda_{i'}(\mathcal{D}_{\tau})+  \delta_n < \epsilon$}} f(\lambda_{i'}(\mathcal{D}_{\tau}) + \delta_n ) \mathcal{U}_{i'} \star \mathcal{U}^{H}_{i'} + \nonumber \\
&  & \sum\limits_{j', \mbox{s.t. $\lambda_{j'}(\mathcal{D}_{\tau}) +  \delta_n 
\geq \epsilon$}} f(\lambda_{j'}(\mathcal{D}_{\tau}) + \delta_n ) \mathcal{U}_{j'} \star \mathcal{U}^{H}_{j'}  \nonumber \\
&\leq &  \sum\limits_{i', \mbox{s.t. $\lambda_{i'}(\mathcal{D}_{\tau}) +  \delta_n < \epsilon$}} f(\lambda_{i'}(\mathcal{D}_{\tau}) ) \mathcal{U}_{i'} \star\mathcal{U}^{H}_{i'} + \nonumber \\
&  & \sum\limits_{j', \mbox{s.t. $\lambda_{j'}(\mathcal{D}_{\tau})  +  \delta_n 
\geq \epsilon$}} f(\lambda_{j'}(\mathcal{D}_{\tau}) + \delta_n ) \mathcal{U}_{j'} \star \mathcal{U}^{H}_{j'}  \nonumber \\
&\leq & f(\mathcal{D}_{\tau}) + \sum\limits_{j', \mbox{s.t. $\lambda_{j'}(\mathcal{D}_{\tau}) +  \delta_n \geq \epsilon$}} f(\lambda_{j'}(\mathcal{D}_{\tau}) + \delta_n ) \mathcal{U}_{j'} \star \mathcal{U}^{H}_{j'}.
\end{eqnarray}
Therefore, the claim in Eq.~\eqref{eq:77} follows by the triangle inequality for $\left\Vert \cdot \right\Vert_{\rho}$ and $f(\lambda_{j'}(\mathcal{D}_{\tau}) + \delta_n )  < \infty$ for $\lambda_{j'}(\mathcal{D}_{\tau}) +  \delta_n \geq \epsilon$. 

\textbf{Eq.~\eqref{eq1:thm:int average thm 10} $\Longleftarrow$ Eq.~\eqref{eq2:thm:int average thm 10}}

The weak majorization relation 
\begin{eqnarray}\label{eq:82}
\prod\limits_{i=1}^{k} \lambda_i (\mathcal{C}) \leq \prod\limits_{i=1}^{k} \exp \int_{\Omega} \log \lambda_i (\mathcal{D}_\tau) d \nu (\tau),
\end{eqnarray}
is valid for $k < m \times p$ from Eq.~\eqref{eq1:thm:weak int average thm 7} $\Longrightarrow$ Eq.~\eqref{eq2:thm:weak int average thm 7} in Theorem~\ref{thm:weak int average thm 7}.  We wish to prove that Eq.~\eqref{eq:82} becomes equal for $k =  m \times p$. It is equivalent to prove that 
\begin{eqnarray}\label{eq:83}
\log \det (\mathcal{C}) \geq   \int_{\Omega} \log \det (\mathcal{D}_{\tau}) d \nu (\tau),
\end{eqnarray}
where $\det ( \cdot )$ is defined by Eq.~\eqref{eq:def T prod tensor determinant}. We can assume that $ \int_{\Omega} \log \det (\mathcal{D}_{\tau}) d \nu (\tau) \geq - \infty$ since Eq.~\eqref{eq:83} is true for  $ \int_{\Omega} \log \det (\mathcal{D}_{\tau}) d \nu (\tau) = - \infty$. Then, $\mathcal{D}_\tau$ are TPD tensors. 

If we scale tensors $\mathcal{C}, \mathcal{D}_{\tau}$ as $a \mathcal{C}, a\mathcal{D}_{\tau}$ by some $a >0$, we can assume $\mathcal{D}_{\tau} \leq \mathcal{I}$ and $\lambda_i(\mathcal{D}_\tau) \leq 1$ for all $ \tau \in \Omega$ and $ i \in \{1,2,\cdots, m \times p\}$. Then for any $p >0$, we have 
\begin{eqnarray}
\frac{1}{m \times p} \left\Vert \mathcal{D}_{\tau}^{-\varrho} \right\Vert_1 \leq \lambda^{-\varrho}_r (\mathcal{D}_{\tau} ) \leq ( \det  (\mathcal{D}_{\tau})  )^{-\varrho},
\end{eqnarray}
and 
\begin{eqnarray}\label{eq:85}
\frac{1}{\varrho} \log \left( \frac{\left\Vert \mathcal{D}^{-\varrho}_\tau \right\Vert_1 }{m \times p}\right)
\leq - \log \det  (\mathcal{D}_{\tau}). 
\end{eqnarray}
If we use tensor trace norm, represented by $\left\Vert \cdot \right\Vert_1$,  as unitarily invariant tensor norm and $f(x) = x^{-\varrho}$ for any $\varrho > 0$ in Eq.~\eqref{eq2:thm:int average thm 10}, we obtain
\begin{eqnarray}\label{eq:86}
\log \left\Vert \mathcal{C}^{-\varrho} \right\Vert_1 \leq \int_{\Omega} \log \left\Vert \mathcal{D}^{-\varrho}_{\tau} \right\Vert_1 d \nu(\tau).
\end{eqnarray}
By adding $\log \frac{1}{m \times p}$ and multiplying $\frac{1}{\varrho}$ for both sides of Eq.~\eqref{eq:86}, we have 
\begin{eqnarray}\label{eq:87}
\frac{1}{\varrho}\log \left( \frac{\left\Vert \mathcal{C}^{-\varrho} \right\Vert_1 }{m \times p} \right)
\leq \int_{\Omega}\frac{1}{\varrho} \log \left(  \frac{\left\Vert \mathcal{D}_\tau^{-\varrho} \right\Vert_1 }{m \times p}           \right) d \nu (\tau)
\end{eqnarray}
Similar to Eqs.~\eqref{eq:52} and~\eqref{eq:53}, we have following two relations as $\varrho \rightarrow 0$:
\begin{eqnarray}\label{eq:88}
\frac{1}{\varrho}\log \left( \frac{\left\Vert \mathcal{C}^{-\varrho} \right\Vert_1 }{m \times p} \right) \rightarrow \frac{- 1}{m \times p} \log \det (\mathcal{C}),
\end{eqnarray}
and 
\begin{eqnarray}\label{eq:89}
\frac{1}{\varrho}\log \left( \frac{\left\Vert \mathcal{D}_{\tau}^{-\varrho} \right\Vert_1}{m \times p}  \right) \rightarrow \frac{- 1}{m \times p} \log \det (\mathcal{D}_\tau).
\end{eqnarray}
From Eq.~\eqref{eq:85} and Lebesgue's dominated convergence theorem, we have 
\begin{eqnarray}\label{eq:90}
\lim\limits_{\varrho  \rightarrow 0} \int_{\Omega}\frac{1}{\varrho}\log \left( \frac{\left\Vert \mathcal{D}_{\tau}^{-\varrho} \right\Vert_1}{m \times p}  \right) d \nu (\tau)= \frac{-1}{m \times p} \int_{\Omega} \log \det(\mathcal{D}_{\tau})     \nu (\tau) 
\end{eqnarray}
Finally, we have Eq.~\eqref{eq:83} from Eqs.~\eqref{eq:87} and~\eqref{eq:90}.

\textbf{Eq.~\eqref{eq1:thm:int average thm 10} $\Longrightarrow$ Eq.~\eqref{eq3:thm:int average thm 10}}

First, we assume that $\mathcal{C}, \mathcal{D}_\tau$ are TPD tensors and $\mathcal{D}_{\tau} \geq \delta \mathcal{I}$ for $\tau \in \Omega$. From Eq.~\eqref{eq1:thm:int average thm 10}, we can apply Theorem~\ref{thm:weak int average thm 5} to $\log \mathcal{C}, \log \mathcal{D}_\tau$ and $f(x) = g(e^x)$ to obtain Eq.~\eqref{eq3:thm:int average thm 10}. 

For $\mathcal{C}, \mathcal{D}_\tau$ are TPSD T-product tensors, we can choose  $\mathbf{a}^{(n)}$ and corresponding  $\mathcal{C}^{(n)}$ for $n \geq n_0$ given $\delta_n \rightarrow 0$ with $\delta_n > 0$ as the proof in Eq.~\eqref{eq1:thm:int average thm 10} $\Longrightarrow$ Eq.~\eqref{eq2:thm:int average thm 10}. Since tensors $\mathcal{C}^{(n)}, \mathcal{D}_\tau + \delta_n \mathcal{I}$ are TPD T-product tensors, we then have 
\begin{eqnarray}\label{eq:92}
\left \Vert g ( \mathcal{C}^{(n)}  ) \right\Vert_{\rho} \leq \int_{\Omega} \left\Vert g ( \mathcal{D}_\tau + \delta_n \mathcal{I} ) \right\Vert_{\rho} d \nu (\tau).
\end{eqnarray}
If $g(0^+) < \infty$, Eq.~\eqref{eq3:thm:int average thm 10} is obtained from Eq.~\eqref{eq:92} by taking $n \rightarrow \infty$. On the other hand, if $g(0^+) = \infty$, we can apply the argument similar to the portion about $f(0^+) = \infty$ in the proof for Eq.~\eqref{eq1:thm:int average thm 10} $\Longrightarrow$ Eq.~\eqref{eq2:thm:int average thm 10} to get $a, b > 0$ such that 
\begin{eqnarray}\label{eq:92 infty}
a \leq \left\Vert g ( \mathcal{D}_\tau + \delta_n \mathcal{I} ) \right\Vert_\rho \leq   \left\Vert  g ( \mathcal{D}_\tau ) \right\Vert_\rho + b,
\end{eqnarray}
for all $\tau \in \Omega$ and $n \geq n_0$. Since the case that $\int_{\Omega} \left\Vert  g ( \mathcal{D}_\tau ) \right\Vert_{\rho} d \nu (\tau) = \infty$ will have Eq.~\eqref{eq3:thm:int average thm 10}, we only consider the case that $\int_{\Omega} \left\Vert  g ( \mathcal{D}_\tau ) \right\Vert_{\rho} d \nu (\tau) < \infty$. Then, we have Eq.~\eqref{eq3:thm:int average thm 10} from Eqs.~\eqref{eq:92},~\eqref{eq:92 infty} and Lebesgue's dominated convergence theorem.

\textbf{Eq.~\eqref{eq1:thm:int average thm 10} $\Longleftarrow$ Eq.~\eqref{eq3:thm:int average thm 10}}

The weak majorization relation
\begin{eqnarray}\label{eq:94}
\sum\limits_{i=1}^{k} \log \lambda_i (\mathcal{C}) \leq 
\sum\limits_{i=1}^{k} \int_{\Omega} \log \lambda_i (\mathcal{D}_\tau) d \nu (\tau)
\end{eqnarray}
is true from the implication from Eq.~\eqref{eq1:thm:weak int average thm 7} to Eq.~\eqref{eq3:thm:weak int average thm 7} in Theorem~\ref{thm:weak int average thm 7}. We have to show that this relation becomes identity for $k=m \times p$. If we apply $\left\Vert \cdot  \right\Vert_{\rho} = \left\Vert \cdot \right\Vert_1$ and $g(x) = x^{-\varrho}$ for any $\varrho > 0$ in Eq.~\eqref{eq3:thm:int average thm 10}, we have 
\begin{eqnarray}\label{eq:95}
\frac{1}{\varrho} \log \left( \frac{ \left\Vert \mathcal{C}^{-\varrho}\right\Vert_1  }{m \times p} \right)
\leq  \frac{1}{\varrho} \log \left( \int_{\Omega} \frac{\left\Vert \mathcal{D}_\tau^{-\varrho}\right\Vert_1}{m \times p}   d \nu(\tau)  \right).
\end{eqnarray}
Then, we will get 
\begin{eqnarray}\label{eq:96}
\frac{- \log \det (\mathcal{C})}{m \times p} &=& \lim\limits_{\varrho \rightarrow 0} \frac{1}{\varrho} \log \left( \frac{ \left\Vert \mathcal{C}^{-\varrho}\right\Vert_1  }{m \times p} \right) \nonumber \\
& \leq & \lim\limits_{\varrho \rightarrow 0} \frac{1}{p} \log \left( \int_{\Omega} \frac{\left\Vert \mathcal{D}_\tau^{-\varrho}\right\Vert_1}{m \times p}   d \nu(\tau)  \right)=_1 \frac{ - \int_{\Omega} \log \det (\mathcal{D}_{\tau}) d \nu (\tau)  }{m \times p},
\end{eqnarray}
which will prove the identity for Eq.~\eqref{eq:94} when $k = m \times p$. The equality in $=_1$ will be proved by the following Lemma~\ref{lma:15}.
$\hfill \Box$

\begin{lemma}\label{lma:15}
Let $\mathcal{D}_\tau$ be TPSD T-product tensors with $\int_{\Omega} \left\Vert \mathcal{D}_{\tau}^{-p}\right\Vert_\rho d \nu (\tau) < \infty$ for any $p > 0$, then we have
\begin{eqnarray}\label{eq1:lma:15}
\lim\limits_{p \rightarrow 0} \left( \frac{1}{p} \log \int_{\Omega} \frac{ \left\Vert \mathcal{D}_\tau^{-p}\right\Vert_1 }{m \times p} d \nu (\tau)\right) &=& -\frac{1}{m \times p} \int_{\Omega} \log \det( \mathcal{D}_{\tau} ) d \nu (\tau)
\end{eqnarray}
\end{lemma}
\textbf{Proof:}
Because $\int_{\Omega} \left\Vert \mathcal{D}_{\tau}^{-p}\right\Vert_\rho d \nu (\tau) < \infty$, we have that $\mathcal{D}_{\tau}$ are TPD tensors for $\tau$ almost everywhere in $\Omega$. Then, we have 
\begin{eqnarray}
\lim\limits_{p \rightarrow 0} \left( \frac{1}{p}\log \int_{\Omega} \frac{ \left\Vert \mathcal{D}_\tau^{-p}\right\Vert_1}{m \times p} d \nu (\tau) \right) &=_1& \lim\limits_{p \rightarrow 0}\frac{   \int_{\Omega} \frac{ - \sum\limits_{i=1}^{m \times p}  \log \lambda_i(\mathcal{D}_\tau)   }{m \times p} d \nu (\tau)   }{  \int_{\Omega} \frac{ \left\Vert \mathcal{D}_\tau^{-p}\right\Vert_1}{m \times p} d \nu (\tau)      } \nonumber \\
&=& \frac{-1}{m \times p} \int_{\Omega} \sum\limits_{i=1}^{m \times p} \log \lambda_i(\mathcal{D}_\tau)  d \nu (\tau)  
\nonumber \\
&=_2& \frac{-1}{m \times p} \int_{\Omega} \log \det ( \mathcal{D}_\tau ) d \nu (\tau), 
\end{eqnarray}
where $=_1$ is from L'Hopital's rule, and $=_2$ is obtained from $\det$ definition.
$\hfill \Box$

\subsection{T-product Tensor Norm Inequalities by Majorization}\label{sec:T-product Tensor Norm Inequalities by Majorization}


In this section, we will apply derived majorization inequalities for T-product tensors to multivariate T-product tensor norm inequalities which will be used to bound random T-product tensor concentration inequalities in later sections. We will begin to present a Lie-Trotter product formula for tensors. 
\begin{lemma}\label{lma: Lie product formula for tensors}
Let $m \in \mathbb{N}$ and $(\mathcal{L}_k)_{k=1}^{m}$ be a finite sequence of bounded T-product tensors with dimensions $\mathcal{L}_k \in  \mathbb{R}^{m \times m \times p}$, then we have
\begin{eqnarray}
\lim_{n \rightarrow \infty} \left(  \prod_{k=1}^{m} \exp(\frac{  \mathcal{L}_k}{n})\right)^{n}
&=& \exp \left( \sum_{k=1}^{m}  \mathcal{L}_k \right)
\end{eqnarray}
\end{lemma}
\textbf{Proof:}

We will prove the case for $m=2$, and the general value of $m$ can be obtained by mathematical induction. 
Let $\mathcal{L}_1, \mathcal{L}_2$ be bounded tensors act on some Hilbert space. Define $\mathcal{C} \define \exp( (\mathcal{L}_1 + \mathcal{L}_2)/n) $, and $\mathcal{D} \define \exp(\mathcal{L}_1/n) \star \exp(\mathcal{L}_2/n)$. Note we have following estimates for the norm of tensors $\mathcal{C}, \mathcal{D}$: 
\begin{eqnarray}\label{eq0: lma: Lie product formula for tensors}
\left\Vert \mathcal{C} \right\Vert, \left\Vert \mathcal{D} \right\Vert \leq \exp \left( \frac{\left\Vert \mathcal{L}_1 \right\Vert + \left\Vert \mathcal{L}_2 \right\Vert  }{n} \right) =  \left[ \exp \left(  \left\Vert \mathcal{L}_1 \right\Vert + \left\Vert \mathcal{L}_2 \right\Vert  \right) \right]^{1/n}.
\end{eqnarray}

From the Cauchy-Product formula, the tensor $\mathcal{D}$ can be expressed as:
\begin{eqnarray}
\mathcal{D} &=& \exp(\mathcal{L}_1/n) \star \exp(\mathcal{L}_2/n) = \sum_{i = 0}^{\infty} \frac{( \mathcal{L}_1/n)^i}{i !} \star \sum_{j = 0}^{\infty} \frac{( \mathcal{L}_2/n)^j}{j !} \nonumber\\
&=& \sum_{l = 0}^{\infty} n^{-l} \sum_{i=0}^l \frac{\mathcal{L}_1^i}{i!} \star \frac{\mathcal{L}_2^{l-i}}{(l - i)!},
\end{eqnarray}
then we can bound the norm of $\mathcal{C} - \mathcal{D}$ as 
\begin{eqnarray}\label{eq1: lma: Lie product formula for tensors}
\left\Vert \mathcal{C} - \mathcal{D} \right\Vert &=& \left\Vert \sum_{i=0}^{\infty} \frac{\left( [ \mathcal{L}_1 + \mathcal{L}_2]/n \right)^i}{i! }
 - \sum_{l = 0}^{\infty} n^{-l} \sum_{i=0}^l \frac{\mathcal{L}_1^i}{i!} \star \frac{\mathcal{L}_2^{l-i}}{(l - i)!} \right\Vert \nonumber \\
&=& \left\Vert \sum_{i=2}^{\infty} k^{-i} \frac{\left( [ \mathcal{L}_1 + \mathcal{L}_2] \right)^i}{i! }
 - \sum_{m = l}^{\infty} n^{-l} \sum_{i=0}^l \frac{\mathcal{L}_1^i}{i!} \star \frac{\mathcal{L}_2^{l-i}}{(l - i)!} \right\Vert \nonumber \\
& \leq & \frac{1}{k^2}\left[ \exp( \left\Vert \mathcal{L}_1 \right\Vert + \left\Vert \mathcal{L}_2 \right\Vert ) + \sum_{l = 2}^{\infty} n^{-l} \sum_{i=0}^l \frac{\left\Vert \mathcal{L}_1 \right\Vert^i}{i!} \cdot \frac{\left\Vert \mathcal{L}_2 \right\Vert^{l-i}}{(l - i)!} \right] \nonumber \\
& = & \frac{1}{n^2}\left[ \exp \left( \left\Vert \mathcal{L}_1 \right\Vert + \left\Vert \mathcal{L}_2 \right\Vert \right) + \sum_{l = 2}^{\infty} n^{-l} \frac{(  \left\Vert \mathcal{L}_1 \right\Vert + \left\Vert \mathcal{L}_2 \right\Vert )^l}{l!} \right] \nonumber \\
& \leq & \frac{2  \exp \left( \left\Vert \mathcal{L}_1 \right\Vert + \left\Vert \mathcal{L}_2 \right\Vert \right) }{n^2}.
\end{eqnarray}

For the difference between the higher power of $\mathcal{C}$ and $\mathcal{D}$, we can bound them as 
\begin{eqnarray}
\left\Vert \mathcal{C}^n - \mathcal{D}^n \right\Vert &=& \left\Vert \sum_{l=0}^{n-1} \mathcal{C}^m (\mathcal{C} - \mathcal{D})\mathcal{C}^{n-l-1} \right\Vert \nonumber \\
& \leq_1 &  \exp ( \left\Vert \mathcal{L}_1 \right\Vert +  \left\Vert \mathcal{L}_2 \right\Vert) \cdot n \cdot \left\Vert \mathcal{L}_1 - \mathcal{L}_2 \right\Vert,
\end{eqnarray}
where the inequality $\leq_1$ uses the following fact 
\begin{eqnarray}
\left\Vert \mathcal{C} \right\Vert^{l} \left\Vert \mathcal{D} \right\Vert^{n - l - 1} \leq \exp \left( \left\Vert \mathcal{L}_1 \right\Vert +  \left\Vert \mathcal{L}_2 \right\Vert \right)^{\frac{n-1}{n}} \leq 
 \exp\left( \left\Vert \mathcal{L}_1 \right\Vert +  \left\Vert \mathcal{L}_2 \right\Vert \right), 
\end{eqnarray}
based on Eq.~\eqref{eq0: lma: Lie product formula for tensors}. By combining with Eq.~\eqref{eq1: lma: Lie product formula for tensors}, we have the following bound
\begin{eqnarray}
\left\Vert \mathcal{C}^n - \mathcal{D}^n \right\Vert &\leq& \frac{2 \exp \left( 2  \left\Vert \mathcal{L}_1 \right\Vert  +  2  \left\Vert \mathcal{L}_2 \right\Vert \right)}{n}.
\end{eqnarray}
Then this lemma is proved when $n$ goes to infity. $\hfill \Box$

Below, new multivariate norm inequalities for T-product tensors are provided according to previous majorization theorems. 
\begin{theorem}\label{thm:Multivaraite Tensor Norm Inequalities}
Let $\mathcal{C}_i \in \mathbb{R}^{m \times m \times p}$ be TPD tensors, where $1 \leq i \leq n$, $\left\Vert \cdot \right\Vert_{\rho}$ be a unitarily invaraint norm with corresponding gauge function $\rho$. For any continous function $f:(0, \infty) \rightarrow [0, \infty)$ such that $x \rightarrow \log f(e^x)$ is convex on $\mathbb{R}$, we have 
\begin{eqnarray}\label{eq1:thm:Multivaraite Tensor Norm Inequalities}
\left\Vert  f \left( \exp \left( \sum\limits_{i=1}^n \log \mathcal{C}_i\right)   \right)  \right\Vert_{\rho} &\leq& \exp \int_{- \infty}^{\infty} \log \left\Vert f \left( \left\vert \prod\limits_{i=1}^{n}  \mathcal{C}_i^{1 + \iota t} \right\vert\right)\right\Vert_{\rho} \beta_0(t) dt ,
\end{eqnarray}
where $\beta_0(t) = \frac{\pi}{2 (\cosh (\pi t) + 1)}$.

For any continous function $g(0, \infty) \rightarrow [0, \infty)$ such that $x \rightarrow g (e^x)$ is convex on $\mathbb{R}$, we have 
\begin{eqnarray}\label{eq2:thm:Multivaraite Tensor Norm Inequalities}
\left\Vert  g \left( \exp \left( \sum\limits_{i=1}^n \log \mathcal{C}_i\right)   \right)  \right\Vert_{\rho} &\leq& \int_{- \infty}^{\infty} \left\Vert g \left( \left\vert \prod\limits_{i=1}^{n}  \mathcal{C}_i^{1 + \iota t} \right\vert\right)\right\Vert_{\rho} \beta_0(t) dt.
\end{eqnarray}
\end{theorem}
\textbf{Proof:}
From Hirschman interpolation theorem~\cite{sutter2017multivariate} and $\theta \in [0, 1]$, we have 
\begin{eqnarray}\label{eq1:Hirschman interpolation}
\log \left\vert h(\theta) \right\vert \leq \int_{- \infty}^{\infty} \log \left\vert h(\iota t) \right\vert^{1 - \theta} \beta_{1 - \theta}(t) d t + \int_{- \infty}^{\infty} \log \left\vert h(1 +  \iota t) \right\vert^{\theta} \beta_{\theta}(t) d t , 
\end{eqnarray}
where $h(z)$ be uniformly bounded on $S \define \{ z \in \mathbb{C}: 0 \leq \Re(z) \leq 1  \}$ and holomorphic on $S$. The term $ d \beta_{\theta}(t) $ is defined as :
\begin{eqnarray}\label{eq:beta theta t def}
\beta_{\theta}(t) \define \frac{ \sin (\pi \theta)}{ 2 \theta (\cos(\pi t) + \cos (\pi \theta))  }.  
\end{eqnarray}
Let $H(z)$ be a uniformly bounded holomorphic function with values in $\mathbb{C}^{m \times m \times p}$. Fix some $\theta \in [0, 1]$ and let $\mathcal{U}, \mathcal{V} \in \mathbb{C}^{m \times m \times p}$ be normalized tensors such that $\langle \mathcal{U}, \mathcal{H}(\theta) \star \mathcal{V} \rangle = \left\Vert H(\theta) \right\Vert$. If we define $h(z)$ as $h(z) \define \langle \mathcal{U}, \mathcal{H}(z) \star \mathcal{V} \rangle $, we have following bound: $\left\vert h(z) \right\vert \leq \left\Vert H(z) \right\Vert $ for all $z \in S$. From Hirschman interpolation theorem, we then have following interpolation theorem for tensor-valued function: 
\begin{eqnarray}\label{eq2:Hirschman interpolation}
\log \left\Vert H(\theta) \right\Vert \leq \int_{- \infty}^{\infty} \log \left\Vert H(\iota t) \right\Vert^{1 - \theta} \beta_{1 - \theta}(t) dt + \int_{- \infty}^{\infty} \log \left\Vert H(1 +  \iota t) \right\Vert^{\theta}  \beta_{\theta}(t) dt .  
\end{eqnarray}

Let $H(z) = \prod\limits_{i=1}^{n} \mathcal{C}^z_i$. Then the first term in the R.H.S. of Eq.~\eqref{eq2:Hirschman interpolation} is zero since $H(\iota t)$ is a product of unitary tensors. Then we have 
\begin{eqnarray}\label{eq:109}
\log \left\Vert \left\vert \prod\limits_{i=1}^{n} \mathcal{C}_i^{\theta} \right\vert^{\frac{1}{\theta}} \right\Vert \leq \int_{- \infty}^{\infty} \log \left\Vert    \prod\limits_{i=1}^{n} \mathcal{C}_i^{1 + \iota t}         \right\Vert \beta_{\theta}(t) d t .  
\end{eqnarray}

From Lemma~\ref{lma:antisymmetric tensor product properties}, we have following relations:
\begin{eqnarray}\label{eq:111-1}
\left\vert \prod\limits_{i=1}^n \left( \wedge^k \mathcal{C}_i \right)^{\theta}\right\vert^{\frac{1}{\theta}} = \wedge^k \left\vert  \prod\limits_{i=1}^n \mathcal{C}^{\theta}_i  \right\vert^{\frac{1}{\theta}},
\end{eqnarray}
and
\begin{eqnarray}\label{eq:111-2}
\left\vert \prod\limits_{i=1}^n \left( \wedge^k \mathcal{C}_i \right)^{1 + \iota t} \right\vert = \wedge^k \left\vert  \prod\limits_{i=1}^n \mathcal{C}^{1 + \iota t}_i  \right\vert.
\end{eqnarray}
If Eq.~\eqref{eq:109} is applied to $\wedge^k \mathcal{C}_i$ for $1 \leq k \leq r$, we have following log-majorization relation from Eqs.~\eqref{eq:111-1} and~\eqref{eq:111-2}:
\begin{eqnarray}\label{eq:112}
\log \vec{\lambda} \left(  \left\vert \prod\limits_{i=1}^{n} \mathcal{C}_i^{\theta} \right\vert^{\frac{1}{\theta}}   \right) \prec \int_{- \infty}^{\infty} \log \vec{\lambda} \left\vert \prod\limits_{i=1}^{n} \mathcal{C}_i^{1 + \iota t } \right\vert^{\frac{1}{\theta}}   \beta_{\theta}(t) d t.
\end{eqnarray}
Moreover, we have the equality condition in Eq.~\eqref{eq:112} for $k = r$ due to following identies:
\begin{eqnarray}\label{eq:113}
\det \left\vert \prod\limits_{i=1}^n \mathcal{C}_i^{\theta} \right\vert^{\frac{1}{\theta}}
= \det \left\vert \prod\limits_{i=1}^n \mathcal{C}_i^{1 + \iota t } \right\vert= \prod\limits_{i=1}^n \det \mathcal{C}_i. 
\end{eqnarray}

At this stage, we are ready to apply Theorem~\ref{thm:int log average thm 10}  for the log-majorization provided by Eq.~\eqref{eq:112} to get following facts:
\begin{eqnarray}\label{eq:114}
\left\Vert  f \left( \left\vert \prod\limits_{i=1}^{n} \mathcal{C}_i^{\theta} \right\vert^{\frac{1}{\theta}}  \right)  \right\Vert_{\rho} &\leq& \exp \int_{- \infty}^{\infty} \log \left\Vert f \left( \left\vert \prod\limits_{i=1}^{n}  \mathcal{C}_i^{1 + \iota t} \right\vert\right)\right\Vert_{\rho} \beta_{\theta}(t) d t ,
\end{eqnarray}
and
\begin{eqnarray}\label{eq:115}
\left\Vert  g \left( \left\vert \prod\limits_{i=1}^{n} \mathcal{C}_i^{\theta} \right\vert^{\frac{1}{\theta}}  \right)  \right\Vert_{\rho} &\leq& \int_{- \infty}^{\infty} \left\Vert g \left( \left\vert \prod\limits_{i=1}^{n}  \mathcal{C}_i^{1 + \iota t} \right\vert\right)\right\Vert_{\rho}  \beta_{\theta}(t) d t.
\end{eqnarray}
From Lie product formula for tensors given by Lemma~\ref{lma: Lie product formula for tensors},  we have 
\begin{eqnarray}\label{eq:117}
 \left\vert \prod\limits_{i=1}^{n} \mathcal{C}_i^{\theta} \right\vert^{\frac{1}{\theta}}
\rightarrow \exp \left(  \sum\limits_{i=1}^{n} \log \mathcal{C}_i  \right). 
\end{eqnarray}
By setting $\theta \rightarrow 0$ in Eqs.~\eqref{eq:114},~\eqref{eq:115} and using Lie product formula given by Eq.~\eqref{eq:117},  we will get Eqs.~\eqref{eq1:thm:Multivaraite Tensor Norm Inequalities} and~\eqref{eq2:thm:Multivaraite Tensor Norm Inequalities}. 
$\hfill \Box$

\section{T-product Tensor Expander Chernoff Bound}\label{sec:T-product Tensor Expander Chernoff Bound}

In this section, we will begin with the derivation for the expectation bound of Ky Fan $k$-norm for the product of TPD tensors in Section~\ref{sec:Expectation Estimation for Product of T-product Tensors}. This bound will be used in the next Section~\ref{sec:Derivation of T-product Tensor Expander Chernoff Bound} for preparing T-product tensor expander Chernoff bounds

\subsection{Expectation Estimation for Product of T-product Tensors}\label{sec:Expectation Estimation for Product of T-product Tensors}

In this section, we will generalize technqiues used in scalar valued expander Chernoff bound proof in~\cite{healy2008randomness} and matrix valued expander Chernoff bound proof in~\cite{garg2018matrix} to build the inequality for the expectation of Ky Fan $k$-norm for the sequential multiplications of TPD T-product tensors. Note that our proof can remove the restriction that the summation of all mapped T-product tensors should be zero tensor, i.e., $\sum\limits_{v \in \mathfrak{V}} \mathrm{g}(v) = \mathcal{O}$.






Let $\mathbf{A}$ be the normalized adjacency matrix of the underlying graph $\mathfrak{G}$ and let $\tilde{\mathbf{A}} = \mathbf{A} \otimes \mathcal{I}_{m^2 \times m^2 \times p^2}$, where the identity tensor $\mathcal{I}_{m^2 \times m^2 \times p^2}$ has dimensions as $m^2 \times m^2 \times p^2$. We use $\mathcal{F} \in \mathbb{R}^{\left(n \times m^2 \right)   \times \left(n  \times m^2 \right) \times p^2  }$ to represent block diagonal T-product tensor valued matrix where the $v$-th diagonal block is the T-product tensor
\begin{eqnarray}
\mathcal{T}_v = \exp\left( \frac{t \mathrm{g}(v) (a + \iota b)}{2} \right) \otimes 
 \exp\left( \frac{t \mathrm{g}(v) (a - \iota b)}{2} \right) \in \mathbb{R}^{m^2   \times  m^2  \times p^2  }. 
\end{eqnarray}
The T-product tensor $\mathcal{F}$ can also be expressed as 
\begin{eqnarray}\label{eq: block decomposition}
\mathcal{F} &=&\left[
    \begin{array}{cccc}
       \mathcal{T}_{v_1}  &  \mathcal{O}  & \cdots &  \mathcal{O}   \\
        \mathcal{O} &  \mathcal{T}_{v_2}  & \cdots & \mathcal{O}  \\
       \vdots &  \vdots & \ddots & \vdots  \\
       \mathcal{O} &  \mathcal{O} & \cdots & \mathcal{T}_{v_n}  \\
    \end{array}
\right].
\end{eqnarray} 
Then the T-product tensor $\left(\mathcal{F} \star  \tilde{\mathbf{A}}\right)^{\kappa}$ is a T-product block tensor valued matrix whose $(u, v)$-block is a tensor with dimensions as $m^2   \times  m^2  \times p^2 $ expressed as :
\begin{eqnarray}
\sum\limits_{v_1, \cdots, v_{\kappa-1} \in \mathfrak{V}} \mathbf{A}_{u,v_1}\left(\prod\limits_{j=1}^{\kappa-2} \mathbf{A}_{v_j, v_{j+1}} \right) 
\mathbf{A}_{v_{\kappa-1}, v} \left(\mathcal{T}_u \star  \mathcal{T}_{v_1} \star \cdots \star \mathcal{T}_{v_{\kappa-1}} \right)
\end{eqnarray}
%

We define a $p^2 \times p^2$ square matrix with first column as all one column, represented by $\mathbb{IO}_{p^2}$. If $p=2$, we have 
\begin{eqnarray}\label{eq:p 2 example}
\mathbb{IO}_{2^2} &=& \left[
    \begin{array}{cccc}
       1  &  0 & 0 &  0  \\
       1  &  0 & 0 &  0  \\
       1  &  0 & 0 &  0  \\
       1  &  0 & 0 &  0  \\
    \end{array}
\right].
\end{eqnarray}
Given a $m \times n$ matrix $\mathbf{X} = [\mathbf{x}_1, \mathbf{x}_2, \cdots, \mathbf{x}_n]$ where $\mathbf{x}_i$ are $m \times 1$ vectors,  we define $\mbox{Vec}(\mathbf{X})$ as:
\begin{eqnarray}
\mbox{Vec}(\mathbf{X}) =  \left[
    \begin{array}{c}
       \mathbf{x}_1 \\
       \mathbf{x}_2 \\
       \vdots  \\
       \mathbf{x}_n  \\
    \end{array}
\right].
\end{eqnarray} 
Then, for T-product tensor $\mathcal{C}, \mathcal{B} \in \mathbb{R}^{m \times m \times p}$, we have following relation:
\begin{eqnarray}\label{eq:kron prod and trace relation}
 \bigl< \mathbb{IO}^\mathrm{T}_{p^2} \otimes   \mbox{Vec}^\mathrm{T}(\mathbf{I}_m) ,   
\mbox{bcirc}(\mathcal{B} \otimes \mathcal{C}) \cdot\left( \mathbb{IO}_{p^2} \otimes   \mbox{Vec}(\mathbf{I}_m) \right)
 \bigr> = \mathrm{Tr}\left( \mathcal{C} \star \mathcal{B}^\mathrm{T} \right),
\end{eqnarray}
where $\cdot$ is the standard matrix multiplication.

Let $\mathbf{U}_0 \in \mathbb{R}^{(n \times m^2 \times p^2) \times p^2}$ be a matrix with $(n \times m^2 \times p^2)$ rows and $p^2$ columns obtained by $\frac{\mathbf{1}}{\sqrt{n}}\otimes \left( \mathbb{IO}_{p^2} \otimes   \mbox{Vec}(\mathbf{I}_m) \right)$, where $\mathbf{1}$ is the all ones vector with size $n$. Then, we will have following expectation of $\kappa$ steps transition of symmetric T-product tensors from the vertex $v_1$ to the vertex $v_{\kappa}$, 
\begin{eqnarray}\label{eq:p17 2nd}
\mathbb{E}\left[\mathrm{Tr} \left(  \prod\limits_{i=1}^{\kappa}\exp \left( \frac{t \mathrm{g}(v_i) (a + \iota b)}{2} \right) \star  \prod\limits_{i=\kappa}^{1}\exp \left( \frac{t \mathrm{g}(v_i) (a -  \iota b)}{2} \right)\right) \right] = ~~~~~~~~~~~~~~~~~~~~~ \nonumber \\
  \biggl<  \mathbf{U}^{\mathrm{T}}_0, \mbox{bcirc} \left( \left( \mathcal{F}\star \tilde{\mathbf{A}} \right)^{\kappa} \right) \cdot \mathbf{U}_0  \biggr>.
\end{eqnarray}
If we define $ \mbox{bcirc} \left( \left( \mathcal{F}\star \tilde{\mathbf{A}} \right)^{\kappa} \right) \cdot \mathbf{U}_0$ as $\mathbf{U}_{\kappa}$, the goal of this section is to estimate $ \bigl<  \mathbf{U}^{\mathrm{T}}_0, \mathbf{U}_{\kappa}  \bigr>$. 

The idea is to separate the space of $\mathbf{U}$ as the subspace spanned by the $m^2 \times p^2$ matrices $\mathbf{1} \otimes e_i$ which is denoted by $\mathbf{U}^{\parallel}$, where $1 \leq i \leq (m^2 \times p^2)$ and $e_i \in \mathbb{R}^{ (m^2 \times p^2) \times p^2 }$ is the matrix with 1 in at the position $i$ of the frst column and 0 elsewhere, and its orthogonal complement space, denoted by $\mathbf{U}^{\perp}$. Following lemma is required to bound how the tensor norm is changed in terms of aforementioned subspace and its orthogonal space after acting by the T-product tensor $\mathcal{F}\star \tilde{\mathbf{A}}$. 

\begin{lemma}\label{lma:4_4}
Given parameters $\lambda \in (0, 1)$, $a  \geq 0$, $r > 0$, and $t > 0$. Let $\mathfrak{G} = (\mathfrak{V}, \mathfrak{E})$ be a regular $\lambda$-expander graph on the vetices set $\mathfrak{V}$ and $\left\Vert \mathrm{g}(v_i) \right\Vert \leq r$ for all $v_i \in \mathfrak{V}$. Each vertex $v \in \mathfrak{V}$ will be assigned a tensor $\acute{\mathcal{T}}_v$, where $\acute{\mathcal{T}}_v \define \frac{ \mathrm{g}(v) (a + \iota b)}{2} \otimes \mathcal{I}_{m,m,p}+ 
 \mathcal{I}_{m,m,p} \otimes \frac{ \mathrm{g}(v) (a - \iota b)}{2} \in \mathbb{R}^{m^2 \times m^2 \times p^2}$. Let $\mathcal{F} \in \mathbb{R}^{\left(n \times m^2 \right)   \times \left(n \times m^2 \right) \times p^2}$ to represent block diagonal T-product tensor valued matrix where the $v$-th diagonal block is the T-product tensor $\exp(t \acute{\mathcal{T}}_v) = \mathcal{T}_v$. For any matrix $\mathbf{U} \in  \mathbb{R}^{(n \times m^2 \times p^2) \times p^2}$, we have following bounds for the spectral norm: 
\begin{enumerate}
   \item $\left\Vert \left( \mbox{bcirc} \left( \mathcal{F} \star\tilde{\mathbf{A}} \right) \cdot \mathbf{U}^{\parallel} \right)^{\parallel} \right\Vert  \leq \gamma_1 \left\Vert  \mathbf{U}^{\parallel} \right\Vert $, where $\gamma_1 = \exp (tr \sqrt{a^2 + b^2}  )$;
   \item $\left\Vert \left( \mbox{bcirc} \left( \mathcal{F} \star\tilde{\mathbf{A}} \right) \cdot \mathbf{U}^{\perp} \right)^{\parallel} \right\Vert  \leq \gamma_2 \left\Vert  \mathbf{U}^{\perp} \right\Vert $, where $\gamma_2 = \lambda (\exp(tr \sqrt{a^2 + b^2}) - 1)$;
   \item $\left\Vert \left( \mbox{bcirc} \left( \mathcal{F} \star\tilde{\mathbf{A}} \right) \cdot \mathbf{U}^{\parallel} \right)^{\perp} \right\Vert  \leq \gamma_3 \left\Vert  \mathbf{U}^{\parallel} \right\Vert $, where $\gamma_3 = \exp(tr \sqrt{a^2 + b^2})-1$;
   \item $\left\Vert \left(  \mbox{bcirc} \left( \mathcal{F} \star\tilde{\mathbf{A}} \right) \cdot \mathbf{U}^{\perp} \right)^{\perp} \right\Vert  \leq \gamma_4 \left\Vert  \mathbf{U}^{\perp} \right\Vert $, where $\gamma_4 = \lambda \exp(tr \sqrt{a^2 + b^2})$.
\end{enumerate}
\end{lemma}
\textbf{Proof:}

For Part 1, let $\mathbf{1} \in \mathbb{R}^n$ be all ones vector, and let $\mathbf{U}^{\parallel} =  \mathbf{1} \otimes \mathbf{v}$ for some $\mathbf{v} \in \mathbb{R}^{ (m^2 \times p^2) \times p^2 }$. Then, we have
\begin{eqnarray}
\left( \mbox{bcirc} \left( \mathcal{F} \star\tilde{\mathbf{A}} \right) \cdot  \mathbf{U}^{\parallel} \right)^{\parallel} 
= \left(\mbox{bcirc} \left( \mathcal{F}  \right) \cdot \mathbf{U}^{\parallel} \right)^{\parallel}
= \mathbf{1} \otimes \left(\frac{1}{n} \sum\limits_{v \in \mathfrak{V}} \mbox{bcirc} \left( \exp(t \acute{\mathcal{T}}_v) \right) \cdot \mathbf{v}  \right)
\end{eqnarray}
and we can bound $ \frac{1}{n} \sum\limits_{v \in \mathfrak{V}} \mbox{bcirc} \left( \exp(t \acute{\mathcal{T}}_v) \right) $ further as 
\begin{eqnarray}
\left\Vert \frac{1}{n} \sum\limits_{v \in \mathfrak{V}} \mbox{bcirc} \left(  \exp(t \acute{\mathcal{T}}_v)  \right) \right\Vert  &=& \left\Vert \frac{1}{n}
\sum\limits_{v \in \mathfrak{V}} \sum\limits_{i=0}^{\infty}  \mbox{bcirc} \left(  \frac{t^i \acute{\mathcal{T}}^i_v }{i !} \right)  \right\Vert  \nonumber \\
&=& \left\Vert \mathbf{I} + \frac{1}{n}\sum\limits_{v \in \mathfrak{V}} \sum\limits_{i=1}^{\infty} 
\mbox{bcirc} \left(  \frac{t^i \acute{\mathcal{T}}^i_v }{i !} \right) \right\Vert  \nonumber \\
&\leq&1 + \frac{1}{n} \sum\limits_{v \in \mathfrak{V}} \sum\limits_{i=1}^{\infty} 
\frac{t^i \left\Vert \mbox{bcirc} \left( \acute{\mathcal{T}}_v \right) \right\Vert^i  }{i !} \nonumber \\
&\leq & 1 + \sum\limits_{i=1}^{\infty}\frac{(tr \sqrt{a^2 + b^2})^i }{i !}=\exp (tr \sqrt{a^2 + b^2}  ),
\end{eqnarray}
where the last inequality is due to the fact that $\left\Vert \frac{t \mathrm{g}(v) (a + \iota b)}{2} \otimes \mathcal{I}_{m,m,p}+ 
 \mathcal{I}_{m,m,p} \otimes \frac{ t \mathrm{g}(v) (a - \iota b)}{2} \right\Vert  \leq 2 tr \times \sqrt{\frac{a^2 + b^2}{4}}  $. 

Then Part 1. of this lemma is established due to 
\begin{eqnarray}
\left\Vert \left( \mbox{bcirc} \left( \mathcal{F} \star\tilde{\mathbf{A}} \right) \cdot \mathbf{U}^{\parallel} \right)^{\parallel} \right\Vert & =& \sqrt{n} \left\Vert \frac{1}{n} \sum\limits_{v \in \mathfrak{V}} \mbox{bcirc} \left( \exp(t \acute{\mathcal{T}}_v) \right) \cdot \mathbf{v}  \right\Vert  \nonumber \\
&\leq & \sqrt{n} \left\Vert \mathbf{v}  \right\Vert  \exp (tr \sqrt{a^2 + b^2}  ) =\exp (tr \sqrt{a^2 + b^2}  )  \left\Vert  \mathbf{U}^{\parallel} \right\Vert .
\end{eqnarray}

For Part 2, since $(\mbox{bcirc} (\tilde{\mathbf{A}}) \cdot \mathbf{U}^{\perp})^{\parallel} = 0$, we have 
\begin{eqnarray}
\left\Vert \left(   \mbox{bcirc} \left( \mathcal{F} \star\tilde{\mathbf{A}} \right) \cdot  \mathbf{U}^{\perp} \right)^{\parallel} \right\Vert  &=& \left\Vert  ( \mbox{bcirc} \left(  ( \mathcal{F} - \mathcal{I}) \star\tilde{\mathbf{A}} \right) \cdot \mathbf{U}^{\perp} )^{\parallel} \right\Vert \nonumber \\
&\leq&   \left\Vert   \mbox{bcirc} \left(  ( \mathcal{F} - \mathcal{I}) \star\tilde{\mathbf{A}} \right) \cdot \mathbf{U}^{\perp} \right\Vert  \nonumber \\
&\leq&  \max\limits_{v \in \mathfrak{V}} \left\Vert  \mbox{bcirc} \left(   \exp(t \acute{\mathcal{T}}_v) - \mathcal{I} \right) \right\Vert   \left\Vert  \mbox{bcirc} \left(   \tilde{\mathbf{A}} \right) \cdot \mathbf{U}^{\perp}   \right\Vert 
\nonumber \\
&\leq&  \max\limits_{v \in \mathfrak{V}} \left\Vert \sum\limits_{i=1}^{\infty} \mbox{bcirc} \left(  \frac{t^i \acute{\mathcal{T}}^i_v }{i !} \right) \right\Vert  \left\Vert  \mbox{bcirc} \left(   \tilde{\mathbf{A}} \right) \cdot \mathbf{U}^{\perp}   \right\Vert \nonumber \\
&\leq&  (\exp (tr \sqrt{a^2 + b^2}  ) -1) \lambda \left\Vert  \mathbf{U}^{\perp} \right\Vert ,
\end{eqnarray}
where the last inequality uses that the underlying graph $\mathfrak{G}$ is a $\lambda$-expander  graph, i.e., $\left\Vert \mathbf{A} \mathbf{x}\right\Vert  \leq \lambda \cdot \left\Vert \mathbf{x}\right\Vert $. Therefore, Part 2 is also valid. 

For Part 3,  because $(\mathbf{U}^{\parallel})^{\perp} = 0$, we have $ \left( \mbox{bcirc} \left(   \mathcal{F} \star\tilde{\mathbf{A}}\right) \cdot \mathbf{U}^{\parallel} \right)^{\perp} = \left( \mbox{bcirc} \left( \mathcal{F} \right) \cdot \mathbf{U}^{\parallel} \right)^{\perp} = ( \mbox{bcirc} \left( \mathcal{F} - \mathcal{I}  \right) \cdot \mathbf{U}^{\parallel }  )^{\perp} $. Then, we can upper bound as 
\begin{eqnarray}
\left\Vert   ( \mbox{bcirc} \left( \mathcal{F} - \mathcal{I}  \right) \cdot \mathbf{U}^{\parallel }  )^{\perp}  \right\Vert &\leq& \left\Vert   \mbox{bcirc} \left( \mathcal{F} - \mathcal{I} \right) \cdot   \mathbf{U}^{\parallel }  \right\Vert  \nonumber \\
& = & \max\limits_{v \in \mathfrak{V}} \left\Vert  \mbox{bcirc} \left(  \exp(t \acute{\mathcal{T}}_v)  - \mathcal{I} \right)  \right\Vert  \cdot \left\Vert \mathbf{U}^{\parallel }  \right\Vert  \nonumber \\
& \leq & \max\limits_{v \in \mathfrak{V}} \left\Vert  \sum\limits_{i=1}^{\infty}  \mbox{bcirc} \left(  \frac{t^i \acute{\mathcal{T}}^i_v }{i !} \right) \right\Vert \cdot \left\Vert \mathbf{U}^{\parallel}   \right\Vert  \nonumber \\
&\leq & (\exp (tr \sqrt{a^2 + b^2}  ) -1)  \left\Vert  \mathbf{U}^{\parallel} \right\Vert,
\end{eqnarray}
hence, Part 3 is also proved. 

Finally, for Part 4, we have 
\begin{eqnarray}
\left\Vert \left( \mbox{bcirc}\left( \mathcal{F} \star\tilde{\mathbf{A}} \right) \cdot\mathbf{U}^{\perp} \right)^{\perp}\right\Vert &\leq& \left\Vert  \mbox{bcirc}\left( \mathcal{F} \star\tilde{\mathbf{A}} \right) \cdot\mathbf{U}^{\perp} \right\Vert  \nonumber \\
& \leq &  \left\Vert \mathcal{F} \right\Vert  \left\Vert \mbox{bcirc}\left(  \tilde{\mathbf{A}} \right) \cdot \mathbf{U}^{\perp} \right\Vert  \leq  \exp (tr \sqrt{a^2 + b^2}  ) \lambda  \left\Vert  \mathbf{U}^{\perp} \right\Vert,
\end{eqnarray}
where we use $ \left\Vert \mathcal{F} \right\Vert \leq  \exp (tr \sqrt{a^2 + b^2}  ) $ (shown at previous part) and the underlying graph $\mathfrak{G}$ is a $\lambda$-expander graph. 
$\hfill \Box$

In the following, we will apply Lemma~\ref{lma:4_4} to bound the following term provided by Eq.~\eqref{eq:p17 2nd}
\begin{eqnarray}
\biggl<  \mathbf{U}^{\mathrm{T}}_0, \mbox{bcirc} \left( \left( \mathcal{F}\star \tilde{\mathbf{A}} \right)^{\kappa} \right) \cdot \mathbf{U}_0  \biggr>
\end{eqnarray}
This bound is formulated by the following Lemma~\ref{lma:4_3}

\begin{lemma}\label{lma:4_3}
Let $\mathfrak{G}$ be a regular $\lambda$-expander graph on the vertex set $\mathfrak{V}$, $g : \mathfrak{V} \rightarrow \mathbb{R}^{m \times m \times p}$, and let $v_1, \cdots, v_{\kappa}$ be a stationary random walk on $\mathfrak{G}$. If $ t r  \sqrt{a^2 + b^2} < 1$ and $\lambda(  2\exp(tr \sqrt{a^2 + b^2}) - 1 ) \leq 1$, we have:
\begin{eqnarray}\label{eq1:lma:4_3}
\mathbb{E}\left[\mathrm{Tr} \left(  \prod\limits_{i=1}^{\kappa}\exp \left( \frac{t \mathrm{g}(v_i) (a + \iota b)}{2} \right) \star  \prod\limits_{i=\kappa}^{1}\exp \left( \frac{t \mathrm{g}(v_i) (a -  \iota b)}{2} \right)\right) \right] \leq  \nonumber \\
(m^2 \times p^2)  \exp \left[ \kappa\left(  2   t r \sqrt{a^2 + b^2} + \frac{8}{1- \lambda} + \frac{16   t r \sqrt{a^2 + b^2} }{1 - \lambda} \right)\right]. ~~~~~~~~~~~~~~~
\end{eqnarray}

\end{lemma}
\textbf{Proof:}
There are two phases for this proof. The first phase is to bound the evolution of tensor norms $\left\Vert \mathbf{U}_i^{\perp}\right\Vert $ and $\left\Vert \mathbf{U}_i^{\parallel}\right\Vert $, respectively. The second phase is to bound $\gamma_i$ for $1 \leq i \leq 4$ in Lemm~\ref{lma:4_4}. We begin with the derivation for the bound $\left\Vert \mathbf{U}_i^{\perp} \right\Vert  $, where $ \mathbf{U}_i$ is the output tensor after acting by the tensor $\mathcal{F} \star\tilde{\mathbf{A}}  $ for $i$ times. It is 
\begin{eqnarray}\label{eq2:lma:4_3}
\left\Vert \mathbf{U}_i^{\perp} \right\Vert  & = & \left\Vert \left(  \mbox{bcirc}\left(\mathcal{F} \star\tilde{\mathbf{A}}  \right) \cdot \mathbf{U}_{i-1} \right)^{\perp} \right\Vert  \nonumber \\
& \leq & \left\Vert \left(   \mbox{bcirc}\left(\mathcal{F} \star\tilde{\mathbf{A}}  \right) \cdot   \mathbf{U}^{\parallel}_{i-1} \right)^{\perp} \right\Vert  + \left\Vert \left(  \mbox{bcirc}\left(\mathcal{F} \star\tilde{\mathbf{A}}  \right) \cdot  \mathbf{U}^{\perp}_{i-1} \right)^{\perp} \right\Vert \nonumber \\
& \leq_1 & \gamma_3 \left\Vert  \mathbf{U}^{\parallel}_{i-1} \right\Vert +
\gamma_4 \left\Vert  \mathbf{U}^{\perp}_{i-1} \right\Vert  \nonumber \\
& \leq_2 & (\gamma_3 + \gamma_3 \gamma_4 +  \gamma_3 \gamma^2_4 + \cdots) \max\limits_{j < i}\left\Vert \mathbf{U}_j^{\parallel} \right\Vert \leq \frac{\gamma_3}{1 - \gamma_4} \max\limits_{j < i}\left\Vert \mathbf{U}_j^{\parallel} \right\Vert,
\end{eqnarray}
where $\leq_1$ is obtained from Lemma~\ref{lma:4_4}, $\leq_2$ is obtained by applying the inequality at $\leq_1$ repeatedly. The next task is to bound  $\left\Vert \mathbf{U}_i^{\parallel} \right\Vert $, we have 
\begin{eqnarray}\label{eq3:lma:4_3}
\left\Vert \mathbf{U}_i^{\parallel} \right\Vert  & = & \left\Vert \left(  \mbox{bcirc}\left(\mathcal{F} \star\tilde{\mathbf{A}}  \right) \cdot    \mathbf{U}_{i-1} \right)^{\parallel} \right\Vert  \nonumber \\
& \leq & \left\Vert  \left(  \mbox{bcirc}\left(\mathcal{F} \star\tilde{\mathbf{A}}  \right) \cdot \mathbf{U}^{\parallel}_{i-1} \right)^{\parallel} \right\Vert  + \left\Vert  \left(  \mbox{bcirc}\left(\mathcal{F} \star\tilde{\mathbf{A}}  \right) \cdot  \mathbf{U}^{\perp}_{i-1} \right)^{\parallel} \right\Vert  \nonumber \\
& \leq_1 & \gamma_1 \left\Vert \mathbf{U}^{\parallel}_{i-1} \right\Vert  +
\gamma_2 \left\Vert  \mathbf{U}^{\perp}_{i-1} \right\Vert  \nonumber \\
& \leq_2 & \left( \gamma_1 + \frac{\gamma_2 \gamma_3}{ 1 - \gamma_4} \right) \max\limits_{j < i}\left\Vert \mathbf{U}_j^{\parallel} \right\Vert,
\end{eqnarray}
where $\leq_1$ is obtained from Lemma~\ref{lma:4_4}, $\leq_2$ is obtained from Eq.~\eqref{eq2:lma:4_3}. From  Eqs~\eqref{eq:p17 2nd},~\eqref{eq2:lma:4_3} and ~\eqref{eq3:lma:4_3}, we have 
\begin{eqnarray}\label{eq4:lma:4_3}
\mathbb{E}\left[\mathrm{Tr} \left(  \prod\limits_{i=1}^{\kappa}\exp \left( \frac{t \mathrm{g}(v_i) (a + \iota b)}{2} \right) \star  \prod\limits_{i=\kappa}^{1}\exp \left( \frac{t \mathrm{g}(v_i) (a -  \iota b)}{2} \right)\right) \right]   \nonumber \\
= \langle \mathbf{U}^{\mathrm{T}}_0, \mathbf{U}_\kappa \rangle 
=  \langle \mathbf{U}^{\mathrm{T}}_0, \mathbf{U}^{\parallel}_\kappa \rangle  \leq  \left\Vert \mathbf{U}_0 \right\Vert \cdot \left\Vert \mathbf{U}^{\parallel}_\kappa \right\Vert = (m \times p) \cdot \left\Vert \mathbf{U}^{\parallel}_\kappa \right\Vert  ~~~~~~~~~~~~~~~ \nonumber \\
\leq   (m \times p) \left( \gamma_1 + \frac{\gamma_2 \gamma_3}{ 1 - \gamma_4} \right)^{\kappa} \cdot \left\Vert \mathbf{U}^{\parallel}_0 \right\Vert 
\leq (m^2 \times p^2) \left( \gamma_1 + \frac{\gamma_2 \gamma_3}{ 1 - \gamma_4} 
\right)^{\kappa}. ~~~~~~~~~~
\end{eqnarray}

The second phase of this proof requires us to bound following four terms: $\gamma_i$ for $1 \leq i \leq 4$. Since $tr \sqrt{a^2 + b^2} < 1$, we can bound $\gamma_1$ as following:
\begin{eqnarray}\label{eq:gamma 1}
\gamma_1 &=& \exp (tr \sqrt{a^2 + b^2}) \leq 1 + 2  t r\sqrt{a^2 + b^2};
\end{eqnarray}
\begin{eqnarray}\label{eq:gamma 2}
\gamma_2 &=& \lambda (\exp (tr \sqrt{a^2 + b^2}) - 1) \leq 2 \lambda   t r\sqrt{a^2 + b^2};
\end{eqnarray}
\begin{eqnarray}\label{eq:gamma 3}
\gamma_3 &=& \exp (tr \sqrt{a^2 + b^2}) - 1 \leq 2   t r \sqrt{a^2 + b^2};
\end{eqnarray}
and the condition $\lambda(  2\exp(tr \sqrt{a^2 + b^2}) - 1 ) \leq 1$, we have 
\begin{eqnarray}\label{eq:gamma 4}
1 - \gamma_4 &  = &1 -  \lambda \exp (t r \sqrt{a^2 + b^2} )  \geq \frac{1 -  \lambda}{2}.
\end{eqnarray}
By applying Eqs.~\eqref{eq:gamma 1},~\eqref{eq:gamma 2},~\eqref{eq:gamma 3} and~\eqref{eq:gamma 4} to the upper bound in Eq.~\eqref{eq4:lma:4_3}, we also have 
\begin{eqnarray}
 (m^2 \times p^2)  \left( \gamma_1 + \frac{\gamma_2 \gamma_3}{ 1 - \gamma_4} 
\right)^{\kappa} \leq  (m^2 \times p^2)  \left[ 1 + 2 (  t r \sqrt{a^2 + b^2} ) + \frac{8 \lambda t^2r^2(a^2 + b^2) }{    1-  \lambda       } \right]^\kappa ~~~~~~~~~~~~~~~  \nonumber \\ 
\leq (m^2 \times p^2)  \left[ \left(1 + 2  t r \sqrt{a^2 + b^2} \right)  \left(1 + \frac{8}{1 - \lambda} \right) \right]^{\kappa} ~~~~~~~~~~~~~~~~~~~~~~  \nonumber \\
\leq (m^2 \times p^2)  \exp \left[ \kappa\left(  2   t r \sqrt{a^2 + b^2} + \frac{8}{1- \lambda} + \frac{16   t r \sqrt{a^2 + b^2} }{1 - \lambda} \right)\right]
\end{eqnarray}
This lemma is proved.
$\hfill \Box$

\subsection{Derivation of T-product Tensor Expander Chernoff Bound}\label{sec:Derivation of T-product Tensor Expander Chernoff Bound}

We begin with a lemma about a Ky Fan $k$-norm inequality for the sum of T-product tensors before deriving our main result about T-product tensor expander Chernoff bound.

\begin{lemma}\label{lma:Ky Fan Inequalities for the sum of tensosrs}
Let $\mathcal{C}_i \in \mathbb{C}^{m \times m \times p}$ be symmetric T-product tensors, then we have 
\begin{eqnarray}\label{eq1:lma:Ky Fan Inequalities for the sum of tensosrs}
\left\Vert \left\vert  \sum\limits_{i=1}^{m} \mathcal{C}_i \right\vert^s \right\Vert_{(k)}
\leq  m^{s -1} \sum\limits_{i=1}^{m}  \left\Vert \left\vert \mathcal{C}_i \right\vert^{s} \right\Vert_{(k)}     
\end{eqnarray}
where $s \geq 1$ and $k \in \{1,2,\cdots, m \times p  \}$. 
\end{lemma}
\textbf{Proof:}
Since we have 
\begin{eqnarray}
\left\Vert \left\vert  \sum\limits_{i=1}^{m} \mathcal{C}_i \right\vert^s \right\Vert_{(k)}
 = \sum\limits_{j=1}^{k} \lambda_j \left( \left\vert  \sum\limits_{i=1}^{m} \mathcal{C}_i \right\vert^s \right) =  \sum\limits_{j=1}^{k} \lambda^s_j \left( \left\vert  \sum\limits_{i=1}^{m} \mathcal{C}_i \right\vert \right) = \sum\limits_{j=1}^{k} \sigma^s_j \left( \sum\limits_{i=1}^{m} \mathcal{C}_i \right). 
\end{eqnarray}
where we have orders for eigenvalues as $\lambda_1 \geq \lambda_2 \geq \cdots $, and singular values as  $\sigma_1 \geq \sigma_2 \geq \cdots $. 

From Lemma 9 in~\cite{sychang2021TProdBernstein} about majorization relation between T-product tensors sum, we have
%
\begin{eqnarray}
\sum\limits_{j = 1}^{k} \sigma_j ( \sum\limits_{i=1}^m \mathcal{C}_i ) \leq 
\sum\limits_{j = 1}^{k} \left( \sum\limits_{i=1}^m \sigma_j (\mathcal{C}_i ) \right),
\end{eqnarray}
where $k \in \{1,2,\cdots, m \times p\}$. Then, we have 
\begin{eqnarray}
\sum\limits_{j = 1}^{k} \sigma^s_j ( \sum\limits_{i=1}^m \mathcal{C}_i ) &\leq &
\sum\limits_{j = 1}^{k} \left( \sum\limits_{i=1}^m \sigma_j (\mathcal{C}_i ) \right)^s
\leq m^{s-1} \sum\limits_{j = 1}^{k} \left( \sum\limits_{i=1}^m \sigma^s_j (\mathcal{C}_i ) \right) \nonumber \\
& = & 
m^{s-1} \sum\limits_{j = 1}^{k} \left( \sum\limits_{i=1}^m \sigma^s_j ( \left\vert \mathcal{C}_i \right\vert ) \right) = m^{s-1} \sum\limits_{j = 1}^{k} \left( \sum\limits_{i=1}^m \sigma_j ( \left\vert \mathcal{C}_i \right\vert^s ) \right) \nonumber \\
& = & m^{s -1} \sum\limits_{i=1}^{m}  \left\Vert \left\vert \mathcal{C}_i \right\vert^{s} \right\Vert_{(k)}    
\end{eqnarray}
$\hfill \Box$

We are ready to present our main theorem about the T-product tensor expander bound for Ky Fan $k$-norm.

\TproductExpanderThm*

\textbf{Proof:}
Let $t > 0$ be a paramter to be chosen later, then we have 
\begin{eqnarray}\label{eq1:thm:tensor expander}
\mathrm{Pr}\left( \left\Vert f \left(  \sum\limits_{j=1}^{\kappa} g(v_j) \right)  \right\Vert_{(k)} \geq \vartheta \right) &=&\mathrm{Pr}\left( \exp\left( \left\Vert t f \left(  \sum\limits_{j=1}^{\kappa} g(v_j) \right)  \right\Vert_{(k)} \right) \geq \exp\left(\vartheta t \right) \right) \nonumber \\
&=_1& \mathrm{Pr}\left(  \left\Vert \exp\left( t f \left(  \sum\limits_{j=1}^{\kappa} g(v_j) \right)  \right)  \right\Vert_{(k)} \geq \exp\left(\vartheta t \right) \right) \nonumber \\
&\leq_2 &   \exp\left( - \vartheta t \right) \mathbb{E} \left(  \left\Vert \exp\left( t f \left(  \sum\limits_{j=1}^{\kappa} g(v_j) \right)  \right)  \right\Vert_{(k)}  \right) \nonumber \\
&\leq_3 &  \exp\left( - \vartheta t \right) \mathbb{E} \left(  \left\Vert f \left( \exp\left( t   \sum\limits_{j=1}^{\kappa} g(v_j) \right)  \right)  \right\Vert_{(k)}  \right) ,
\end{eqnarray}
where equality $=_1$ comes from spectral mapping theorem and TPD of , inequality $\leq_2$ is obtained from Markov inequality, and the last inequality $\leq_3$ is based on our function $f$ assumption (first assumption).

From Eq.~\eqref{eq2:thm:Multivaraite Tensor Norm Inequalities} in Theorem~\ref{thm:Multivaraite Tensor Norm Inequalities}, we can further bound the expectation term in Eq.~\eqref{eq1:thm:tensor expander} as 
\begin{eqnarray}\label{eq2:thm:tensor expander}
\mathbb{E} \left(  \left\Vert f \left( \exp\left( t   \sum\limits_{j=1}^{\kappa} g(v_j) \right)  \right)  \right\Vert_{(k)}  \right)  ~~~~~~~~~~~~~~~~~~~~~~~~~~~~~~~~~~~~~~~~~~~~~~~~~~~~~~~~~~~~~~~~~~~~~~~~~~~~~ \nonumber \\
 \leq \mathbb{E} \left( \int\limits_{- \infty}^{\infty} \left\Vert f \left(   \left\vert \prod\limits_{j=1}^{\kappa} \exp\left( t g(v_j) (1 + \iota \tau   ) \right) \right\vert  \right) \right\Vert_{(k)}\beta_0( \tau) d   \tau   \right) ~~~~~~~~~~~~~~~~~~~~~~~~~~~~~~~~~ \nonumber \\
=_1  \mathbb{E} \left( \int\limits_{- \infty}^{\infty} \left\Vert \left(
\sum\limits_{l=0}^{n}a_l   \left\vert \prod\limits_{j=1}^{\kappa} \exp\left( t g(v_j) (1 + \iota \tau   ) \right) \right\vert^l 
 \right)^s
 \right\Vert_{(k)}  \beta_0( \tau) d \tau  \right)  ~~~~~~~~~~~~~~~~~~~~~~ \nonumber \\
\leq_2 (n+1)^{(s-1)} \mathbb{E} \left( \int\limits_{- \infty}^{\infty} \sum\limits_{l=0}^{n}a_l  \left\Vert \left(
  \left\vert \prod\limits_{j=1}^{\kappa} \exp\left( t g(v_j) (1 + \iota \tau   ) \right) \right\vert^l 
 \right)^s
 \right\Vert_{(k)} \beta_0( \tau) d   \tau   \right)    \nonumber \\
= (n+1)^{(s-1)} \cdot ~~~~~~~~~~~~~~~~~~~~~~~~~~~~~~~~~~~~~~~~~~~~~~~~~~~~~~~~~~~~~~~~~~~~~~~~~~~~~~~~~~~~~~~~~~~~~~~~~~~~~~ \nonumber \\
\left( \int\limits_{- \infty}^{\infty} \sum\limits_{l=0}^{n}a_l  \mathbb{E} \left( \left\Vert \left(
  \left\vert \prod\limits_{j=1}^{\kappa} \exp\left( t g(v_j) (1 + \iota \tau   ) \right) \right\vert^l 
 \right)^s
 \right\Vert_{(k)} \right)  \beta_0( \tau) d   \tau   \right) 
\end{eqnarray}
where equality $=_1$ comes from the function $f$ definition, inequality $\leq_2$ is based on Lemma~\ref{lma:Ky Fan Inequalities for the sum of tensosrs}. Each summand for $l \ge 1$ in Eq.~\eqref{eq2:thm:tensor expander} can further be bounded as 
\begin{eqnarray}\label{eq3:thm:tensor expander}
 \mathbb{E} \left( \left\Vert \left(
  \left\vert \prod\limits_{j=1}^{\kappa} \exp\left( t g(v_j) (1 + \iota \tau   ) \right) \right\vert^l 
 \right)^s
 \right\Vert_{(k)} \right)  = \mathbb{E} \left( \left\Vert 
  \left\vert \prod\limits_{j=1}^{\kappa} \exp\left( t g(v_j) (1 + \iota \tau   ) \right) \right\vert^{ls} 
  \right\Vert_{(k)} \right) \nonumber \\
\leq_1   \frac{1}{mp} \mathbb{E} \left(  \mathrm{Tr}\left(   \left\vert \prod\limits_{j=1}^{\kappa} \exp\left( t g(v_j) (1 + \iota \tau   ) \right) \right\vert^{ls} \right) \right) + ~~~~~~~~~~~~~~~~~~~~~~~~~~~~~~~~~~~~~~~~~~~~~ \nonumber \\
  \mathbb{E} \left(  \sqrt{\frac{(mp - k)}{k mp }   \mathrm{Tr} \left(  \left\vert \prod\limits_{j=1}^{\kappa} \exp\left( t g(v_j) (1 + \iota \tau   ) \right) \right\vert^{2ls} \right)} \right) ~~~~~~~~~~~~~~~~~~~~~~~~~~~~~~~~~~~~~ \nonumber \\
\leq_2 mp 
\exp \left[ \kappa\left(  2   lst r \sqrt{1 + \tau^2} + \frac{8}{1- \lambda} + \frac{16   tls r \sqrt{1 + \tau^2} }{1 - \lambda} \right)\right] + ~~~~~~~~~~~~~~~~~~~~~~ \nonumber \\ 
\left\{ \frac{(mp - k)mp}{k}\exp \left[ \kappa\left(  4 lst r \sqrt{1 + \tau^2} + \frac{8}{1- \lambda} + \frac{32  l s t r \sqrt{1 + \tau^2} }{1 - \lambda} \right)\right] \right\}^{1/2} ~~~~~~ \nonumber \\
\leq_3 \left( mp+ \sqrt{ \frac{(mp-k)mp}{k}} \right)\cdot \exp\left[ \kappa\left(  2   lst r (1 +\tau) + \frac{8}{1- \lambda} + \frac{16   tls r (1 +\tau)  }{1 - \lambda} \right)\right] 
\end{eqnarray}
where $\leq_1$ comes from Theorem 2.2 in~\cite{wolkowicz1980bounds} and our T-product tensor trace definition provided by Eq.~\eqref{eq:trace def}, and $\leq_2$ comes from Lemma~\ref{lma:4_3}, and the last inequality $\leq_3$ is obtaiend by bounding $\sqrt{1 + \tau^2}$ as $1 + \tau$.

From Eqs.~\eqref{eq1:thm:tensor expander},~\eqref{eq2:thm:tensor expander}, and~\eqref{eq3:thm:tensor expander}, we have 
\begin{eqnarray}\label{eq4:thm:tensor expander}
\mathrm{Pr}\left( \left\Vert f \left(  \sum\limits_{j=1}^{\kappa} g(v_j) \right)  \right\Vert_{(k)} \geq \vartheta \right) \leq \min\limits_{t > 0   } \left[ (n+1)^{(s-1)} e^{-\vartheta t} \left(a_0 k  + \left(  mp + \sqrt{ \frac{(mp-k)mp}{k} }  \right)\cdot \right. \right.  \nonumber \\
\left. \left. \sum\limits_{l=1}^n a_l \int\limits_\infty^{\infty} 
\exp\left[ \kappa\left(  2   lst r (1 +\tau) + \frac{8}{1- \lambda} + \frac{16   tls r (1 +\tau)  }{1 - \lambda} \right)\right] \beta_0(\tau) d \tau \right)\right] ~~~~~~~ \nonumber \\
\leq_1  \min\limits_{t > 0   } \left[ (n+1)^{(s-1)} e^{-\vartheta t} \left(a_0 k  + \left( mp + \sqrt{ \frac{(mp-k)mp}{k} }  \right)\cdot \right. \right.  ~~~~~~~~~~~~~~~~~~~~~~~~~~~~~~~~~~~~~~~ \nonumber \\
\left. \left. \sum\limits_{l=1}^n a_l \int\limits_\infty^{\infty} 
\exp\left[ \kappa\left(  2   lst r (1 +\tau) + \frac{8}{1- \lambda} + \frac{16   tls r (1 +\tau)  }{1 - \lambda} \right)\right]  \frac{C \exp( \frac{-\tau^2}{2 \sigma^2} ) }{\sigma \sqrt{2 \pi}} d \tau  \right)\right] \nonumber \\
 =  \min\limits_{t > 0 } \left\{ (n+1)^{(s-1)} e^{-\vartheta t} \left[ a_0 k  +C \left(  mp + \sqrt{ \frac{(mp-k)mp}{k} }  \right)\cdot \right. \right. ~~~~~~~~~~~~~~~~~~~~~~~~~~~~~~~~~~~~~ \nonumber \\
\left. \left. \sum\limits_{l=1}^n a_l \exp\left( 8 \kappa \overline{\lambda} + 2  (\kappa +8 \overline{\lambda}) lsr t + 2(\sigma (\kappa +8 \overline{\lambda}) lsr )^2 t^2 \right)  \right] \right\},~~~~~~~~~~~~~~~~~~~
\end{eqnarray}
where inequality $\leq_1$ is obtained by the distribution bound for $ \beta_0(\tau)$ via another distribution function $\frac{C \exp( \frac{-\tau^2}{2 \sigma^2} ) }{\sigma \sqrt{2 \pi}}$, and the last equality comes from Gaussian integral with respect to the variable $\tau$ by setting $1 -\lambda$ as $\overline{\lambda}$. 
$\hfill \Box$

Following corollary is about a tensor expander bound with identity function $f$. 

\begin{corollary}\label{cor:tensor expander one variable}
If we consider the special case of Theorem~\ref{thm:tensor expander} by assuming that the function $f: x \rightarrow x$ is an identiy map, then we have 
\begin{eqnarray}
\mathrm{Pr}\left( \left\Vert  \sum\limits_{j=1}^{\kappa} g(v_j)  \right\Vert_{(k)} \geq  \vartheta   \right) &\leq &  C \left(  mp + \sqrt{ \frac{(mp-k)mp}{k} } \right) \cdot \nonumber \\
&  & \exp \left( - \frac{ \vartheta ^2}{8 \sigma^2 r^2} + \frac{ \vartheta }{2 \sigma^2 r^2} - \frac{1}{2\sigma^2} + 8\kappa \overline{\lambda} \right).
\end{eqnarray}
\end{corollary}
\textbf{Proof:}
From Theorem~\ref{thm:tensor expander}, since the exponent is a quadratic function of $t$, the minimum of this quadratic function is achieved by selecting $t$ as
\begin{eqnarray}\label{eq1:cor:tensor expander one variable}
t = \frac{\vartheta - 2 (\kappa + 8 \overline{\lambda}) r}{4 \sigma^2 r^2 (\kappa + 8 \overline{\lambda} )^2},
\end{eqnarray}
then, we have the desired bound after some algebra by applying Eq.~\eqref{eq1:cor:tensor expander one variable} in Eq.~\eqref{eq0:thm:tensor expander} and setting $l=s=1$, all $a_i=0$ for $1 \leq i \leq n$ except $a_1 = 1$. 
$\hfill \Box$

\section{Conclusions}\label{sec:Conclusions}

In this work, we first build tensor norm inequalities for T-product tensors based on the concept of log-majorization, and apply these new tensor norm inequalities to derive the \emph{T-product tensor expander Chernoff bound} which generalizes the matrix expander Chernoff bound by adopting more general norm for tensors, Ky Fan norm, and general convex function, instead of identity function, of random T-product tensors summation. 

There are several future directions that can be explored based on the current work. The first is to consider other types of T-product tensor expander Chernoff bound under other non-independent assumptions among random T-product tensors. Another direction is to characterize random behaviors of other T-product tensor related quantities besides norms or eigenvalues, for example, what is the T-product tensor rank behavior for the summation of random T-prduct tensors.


\bibliographystyle{IEEETran}
\bibliography{TProd_Major_Expander_Bib}

\end{document}